\documentclass{amsart}
\usepackage{pst-plot}
\usepackage{amsmath}
\usepackage{amssymb}
\usepackage{mathrsfs}

\newtheorem{theorem}{Theorem}[section]
\newtheorem{corollary}[theorem]{Corollary}
\newtheorem{proposition}[theorem]{Proposition}
\newtheorem{lemma}[theorem]{Lemma}
\theoremstyle{definition}
\newtheorem{definition}[theorem]{Definition}
\newtheorem{example}[theorem]{Example}

\theoremstyle{remark}

\numberwithin{equation}{section}

\newcommand{\qbinom}[2]{\begin{bmatrix} #1 \\ #2 \end{bmatrix}_q}

\title{Overpartitions, lattice paths, and Rogers-Ramanujan identities}

\author{Sylvie Corteel}
\address{LRI, CNRS et Universit\'e Paris-Sud, B\^at 490, F-91405 Orsay France}
\email{corteel@lri.fr}
\thanks{The authors are partially supported by the ACI Jeunes Chercheurs
``Partitions d'entiers \`a la fronti\`ere de la combinatoire, des
$q$-s\'eries et de la th\'eorie des nombres.''}

\author{Olivier Mallet}
\address{LIAFA, Universit\'e Denis Diderot et CNRS, 
2 place Jussieu, Case 7014, F-75251 Paris Cedex 05}
\email{mallet@liafa.jussieu.fr}



\begin{document}

\begin{abstract}
We extend partition-theoretic work of Andrews, Bressoud, and Burge
to overpartitions, defining the notions of successive ranks,
generalized  Durfee squares, and generalized lattice paths, and
then relating these to overpartitions defined by multiplicity
conditions on the parts.  This leads to many new partition and
overpartition identities, and provides a unification of a number
of well-known identities of the Rogers-Ramanujan type.  Among
these are Gordon's generalization of the Rogers-Ramanujan
identities, Andrews' generalization of the G\"ollnitz-Gordon
identities, and Lovejoy's ``Gordon's theorems for overpartitions."
%We show how these combinatorial statistics give extensions to overpartitions of combinatorial
%interpretations in terms of lattice paths of the generalizations of the Rogers-Ramanujan identities
%due to Burge, Andrews and Bressoud.
%All our proofs are combinatorial and use bijective techniques. Our
%results
%include the Andrews-Gordon identities, the generalization of the Gordon-G\"ollnitz identities
%and Gordon's theorems for overpartitions.
\end{abstract}

\maketitle

\section{Introduction}
In 1961 Gordon established his celebrated combinatorial
generalization of the Rogers-Ramanujan identities:
\begin{theorem}\cite{go} \label{Gordonthm}
Let $B_{k,i}(n)$ denote the number of partitions of $n$ of the
form $(\lambda_1,\lambda_2,\ldots ,\lambda_s)$, where
$\lambda_\ell -\lambda_{\ell+k-1}\ge  2$ and at most $i - 1$ of
the parts are equal to 1. Let $A_{k,i}(n)$ denote the number of
partitions of $n$ into parts not congruent to $0, \pm i$ modulo
$2k + 1$. Then $A_{k,i}(n) = B_{k,i}(n)$.
\end{theorem}
%Gordon's theorem is an extension of the famous Rogers-Ramanujan
%identities proved by Rogers as a $q$-series identity \cite{ro} and
%interpreted combinatorially by MacMahon \cite{MM}.  These
%identities correspond to the cases $k=2$ and $i=1$ or $2$.
Over the years, a number of other combinatorial functions have
been found to be equal to the $A_{k,i}(n)$ in Gordon's theorem.
Most notable, perhaps, are two results of Andrews that employ
Atkin's successive ranks \cite{at} and Andrews' own new
idea of Durfee dissection :
\begin{theorem}\cite{ag3} \label{andrews-gordonthm}
Let $C_{k,i}(n)$ be the number of partitions of $n$ whose
successive ranks lie in the interval $[-i+2,2k-i-1]$ and let
$D_{k,i}(n)$ be the number of partitions of $n$ with $i-1$
successive Durfee squares followed by $k-i$ successive Durfee
rectangles. Then
$$
A_{k,i}(n) = B_{k,i}(n)=C_{k,i}(n) = D_{k,i}(n).
$$
\end{theorem}

An overpartition is a partition where the final occurrence of a
part can be overlined \cite{cl}. For example, there exist 8
overpartitions of 3
$$
(3),\ (\overline{3}),\ (2,1),\ (\overline{2},1),\ (2,\overline{1}),\
(\overline{2},\overline{1}),\ (1,1,1),\ (1,1,\overline{1}).
$$
In recent years overpartitions have been heavily studied,
sometimes under different names and guises. They have been called
joint partitions \cite{bp1}, or dotted partitions \cite{br}, and
they are closely related to 2-modular diagrams \cite{mm}, jagged
partitions \cite{fm,fm1} and superpartitions \cite{l}.
Overpartitions arise in the study of the combinatorics of basic
hypergeometric series identities \cite{cl1,fm1,lo,lo1,y},
congruences properties of modular forms \cite{fm,ma},
supersymmetric functions \cite{l}, Lie algebras \cite{kw} and mathematical physics \cite{l,fm,fm1}.

In 2003 Lovejoy \cite{lo} proved an overpartition identity wherein
one of the functions closely resembles the $B_{k,i}(n)$ in
Gordon's theorem:
\begin{theorem}[``Gordon's theorem for overpartitions"] \cite{lo}
\label{Gordonoverthm} Let $\overline{B}_k(n)$ denote the number of
overpartitions of $n$ of the form $(\lambda_1,\lambda_2,\ldots
,\lambda_s)$, where $\lambda_\ell - \lambda_{\ell+k-1} \ge 1$ if
$\lambda_{\ell+k-1}$ is overlined and $\lambda_\ell -
\lambda_{\ell+k-1} \ge 2$ otherwise. Let $\overline{A}_k(n)$
denote the number of overpartitions of $n$ into parts not
divisible by $k$. Then $\overline{A}_k(n) = \overline{B}_k(n)$.
\end{theorem}
Notice that Lovejoy's result can be seen as an overpartition
analogue of Gordon's theorem, in the sense that the conditions on
the $\overline{B}_k(n)$ reduce to the conditions on the
$B_{k,k}(n)$ if the overpartition has no overlined parts and is
indeed a partition.

Two questions naturally arise.  First, given the similarities
between Theorems \ref{Gordonthm} and \ref{Gordonoverthm}, is there
some general framework which contains these two theorems?  Second,
is there an analogue for overpartitions of Andrews' result,
Theorem \ref{andrews-gordonthm}?  In this paper, both of these
questions shall be answered in the affirmative.  Moreover, our
results contain many other partition and overpartition identities,
such as Andrews' generalization of the G\"ollnitz-Gordon
identities \cite{ag4}.

%Why do you write this?
%It is
%still a well known open problem to find a natural bijective proof
%of these identities, even though an impressive number of nearly
%combinatorial proofs have been published. A recent example is
%\cite{bp}.

It is well understood combinatorially that $B_{k,i}(n)=C_{k,i}(n)
= D_{k,i}(n)$ and this result was established by some beautiful
work of Burge \cite{bu1,bu2} using some recursive arguments. This
work was reinterpreted by Andrews and Bressoud \cite{ab} who
showed that Burge's argument could be rephrased in terms of lattice paths with 
two kinds of steps and that Gordon's theorem can be established thanks to these
combinatorial arguments and the Jacobi Triple product identity
\cite{gr}. Finally Bressoud \cite{b} reinterpreted these in terms
of different lattice paths with three kinds of steps and 
gave some direct bijections between the
objects counted by $B_{k,i}(n)$, $C_{k,i}(n)$, $D_{k,i}(n)$ and
the lattice paths.

With our main theorem we extend the main results of the above
works \cite{ab,b,bu1,bu2} to overpartitions.  In particular, we
generalize all of the combinatorial definitions to overpartitions
and successfully adapt the methods of proof.  This is the result
that provides a unifying framework for Theorems \ref{Gordonthm},
\ref{andrews-gordonthm}, and \ref{Gordonoverthm}.

\begin{theorem}
~
\begin{itemize}
\item Let $\overline{B}_{k,i}(n,j)$ be the number of overpartitions of $n$ of the form $(\lambda_1,\lambda_2,\ldots ,\lambda_s)$ with $j$ overlined parts, where
$\lambda_{\ell} - \lambda_{\ell+k-1} \ge 1$ if $\lambda_{\ell+k-1}$ is overlined and $\lambda_{\ell} -
\lambda_{\ell+k-1} \ge 2$ otherwise
and at most $i-1$ parts are equal to $1$.
\item Let $\overline{C}_{k,i}(n,j)$ be the number of
overpartitions of $n$ with $j$ non-overlined parts in the bottom row of their Frobenius representation
and whose  successive ranks
lie in $[-i+2,2k-i-1]$.
\item Let $\overline{D}_{k,i}(n,j)$ be the number of overpartitions of $n$ with
$j$ overlined parts and $i-1$ successive Durfee
squares followed by $k-i$ successive Durfee rectangles, the first one being a generalized Durfee square/rectangle.
\item Let $\overline{E}_{k,i}(n,j)$ be the number of  paths that use four kinds of unitary steps with
special $(k,i)$-conditions, major index $n$, and $j$ South steps.
\end{itemize}
Then $\overline{B}_{k,i}(n,j)=\overline{C}_{k,i}(n,j)=\overline{D}_{k,i}(n,j)=\overline{E}_{k,i}(n,j)$.
\label{main}
\end{theorem}

All of the combinatorial notions in this theorem will be defined
in detail in Section 2.  The addition of the generalized lattice
paths counted by $\overline{E}_{k,i}(n,j)$ is the key step which allows us to
prove Theorem \ref{main} combinatorially.  In terms of generating
functions, we have:
\begin{theorem} The generating function $\overline{\mathcal
E}_{k,i}(a,q)=\sum_{n,j}\overline{E}_{k,i}(n,j)q^na^j$ is:
\begin{equation}
\overline{\mathcal E}_{k,i}(a,q)=
\frac{(-aq)_\infty}{(q)_\infty}\sum_{n=-\infty}^{\infty}(-1)^n a^n
q^{kn^2+(k-i+1)n}\frac{(-1/a)_n}{(-aq)_n}. \label{gf}
\end{equation}
\label{thegf}
\end{theorem}

Here we have used the classical $q$-series notations:
\begin{eqnarray*}
(a)_\infty &=&(a;q)_\infty=\prod_{i=0}^{\infty}(1-aq^i),\\
(a)_n&=&(a)_\infty/(aq^n)_\infty,\\
(a_1,\ldots,a_k;q)_\infty &=& (a_1;q)_\infty\ldots (a_k;q)_\infty.
\end{eqnarray*}

In several cases, we can use the Jacobi Triple Product identity
\cite{gr}:
\begin{equation}
(-1/z,-zq,q;q)_\infty=\sum_{n=-\infty}^{\infty} z^n q^{\binom{n+1}{2}}
\label{JTP}
\end{equation}
to show that this generating function has a very nice form. For
example,
\begin{corollary}
\begin{eqnarray}
\overline{\mathcal E}_{k,i}(0,q)&=&\frac{(q^i,q^{2k+1-i},q^{2k+1};q^{2k+1})_\infty}{(q)_\infty}\label{andrews}\\
\overline{\mathcal E}_{k,i}(1/q,q^2)&=&\frac{(q^2;q^4)_\infty(q^{2i-1},q^{4k+1-2i},q^{4k};q^{4k})_\infty}{(q)_\infty}\label{gordon}\\
\overline{\mathcal E}_{k,i}(1,q)&=&\frac{(-q)_\infty}{(q)_\infty}
\sum_{j=0}^{2(k-i)}(-1)^j(q^{i+j},q^{2k-i-j},q^{2k};q^{2k})_\infty
\label{gord}\\
\overline{\mathcal E}_{k,i}(1/q,q)&=&\frac{(-q)_\infty}{(q)_\infty}(
(q^i,q^{2k-i},q^{2k};q^{2k})_\infty\nonumber\\
&&\ \ \ \ \ \ \ \ +(q^{i-1},q^{2k+1-i},q^{2k};q^{2k})_\infty)\label{gord1}
\end{eqnarray}
\label{aftergf}
\end{corollary}

Hence our result gives a general view of different
problems on partitions and overpartitions and shows how they are related.
\begin{itemize}
\item Equation (\ref{andrews}) corresponds to the Andrews-Gordon identities \cite{ag3}.
\item Equation (\ref{gordon}) corresponds to Andrews's generalization of
the Gordon-G\"ollnitz identities \cite{ag4,ab}.
\item Equation (\ref{gord}) with $i=k$ and (\ref{gord1}) with $i=1$
correspond to the two Gordon's theorems for overpartitions of Lovejoy \cite{lo}.
\end{itemize}
Therefore our extension of the work on the Andrews-Gordon identities \cite{ab,b,bu1,bu2}
to the case of overpartitions includes these identities, but it also includes Andrews's generalization of
the Gordon-G\"ollnitz identities and Gordon's theorems for 
overpartitions. We prove this Corollary  and deduce some new partition theorems in Section 7.

In Section 2 we present all the necessary notions.
In Section 3 we present the paths counted by $\overline{E}_{k,i}(n,j)$ and compute the generating function. In Section 4 we present a direct bijection
between the paths counted by $\overline{E}_{k,i}(n,j)$ and the overpartitions counted by $\overline{C}_{k,i}(n,j)$.
In Section 5 we present a recursive bijection
between the paths counted by $\overline{E}_{k,i}(n,j)$ and the overpartitions
counted by $\overline{B}_{k,i}(n,j)$. We also give a generating function proof. In Section 6, we present
a combinatorial argument that shows that the paths counted by $\overline{E}_{k,i}(n,j)$ and the overpartitions
counted by $\overline{D}_{k,i}(n,j)$ are in bijection. All these bijections are refinements of Theorem
\ref{main}. The number of peaks of the paths will correspond respectively
to the number of columns of the Frobenius representations, the length of the multiplicity sequence and the size of the generalized Durfee square. 
In Section 7 we prove Corollary \ref{aftergf} and interpret 
it combinatorially.
We conclude in Section 8
with further questions.

\section{Definitions on overpartitions}

We will define all the notions in terms of overpartitions. We refer to \cite{ag} for definitions for partitions.
In all of the cases the definitions coincide when the overpartition has no overlined parts.

\begin{figure}[ht]
        \centering
        \psset{unit=0.4}
\begin{pspicture}(0,0)(5,4)
% part 5
\rput(0.5,3.5){2}
\rput(1.5,3.5){\small 2}
\rput(2.5,3.5){\small 2}
\rput(3.5,3.5){\small 2}
\rput(4.5,3.5){\small 1}
\psframe[dimen=middle](0,4)(1,3)
\psframe[dimen=middle](1,4)(2,3)
\psframe[dimen=middle](2,4)(3,3)
\psframe[dimen=middle](3,4)(4,3)
\psframe[dimen=middle](4,4)(5,3)
% part 4
\rput(0.5,2.5){\small 2}
\rput(1.5,2.5){\small 2}
\rput(2.5,2.5){\small 2}
\rput(3.5,2.5){\small 2}
\psframe[dimen=middle](0,3)(1,2)
\psframe[dimen=middle](1,3)(2,2)
\psframe[dimen=middle](2,3)(3,2)
\psframe[dimen=middle](3,3)(4,2)
% part 3
\rput(0.5,1.5){\small 2}
\rput(1.5,1.5){\small 2}
\rput(2.5,1.5){\small 2}
\psframe[dimen=middle](0,2)(1,1)
\psframe[dimen=middle](1,2)(2,1)
\psframe[dimen=middle](2,2)(3,1)
% part 3
\rput(0.5,0.5){\small 2}
\rput(1.5,0.5){\small 2}
\rput(2.5,0.5){\small 1}
\psframe[dimen=middle](0,1)(1,0)
\psframe[dimen=middle](1,1)(2,0)
\psframe[dimen=middle](2,1)(3,0)
\end{pspicture}
        \caption{Example of a 2-modular diagram of weight 28.}
        \label{2modular}
\end{figure}
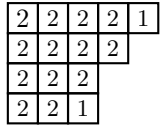

An overpartition of $n$ is a non-increasing sequence of natural
numbers whose sum is $n$ in which the final occurrence
(equivalently, the first occurrence) of a number may be overlined.
Alternatively $n$ can be called the weight of the overpartition.
Since the overlined parts form a partition into distinct
parts and the non-overlined parts form an ordinary partition, the generating function for overpartitions is
$\frac{(-q)_\infty}{(q)_\infty}$.
The Ferrers diagram of an overpartition is a classical Ferrers diagram
where the corners can be marked (see Figure \ref{Ferrers}). A 2-modular diagram is a Ferrers diagram
of an overpartition where the marked corners are filled with ones and
the other cells are filled with twos (see Figure \ref{2modular}). The weight is the sum of the entries.

\begin{figure}[ht]
        \centering
        \psset{unit=0.4}
\begin{pspicture}(0,0)(5,4)
% part 5 (surlignee)
\psframe[dimen=middle](0,4)(1,3)
\psframe[dimen=middle](1,4)(2,3)
\psframe[dimen=middle](2,4)(3,3)
\psframe[dimen=middle](3,4)(4,3)
\psframe[dimen=middle,fillstyle=solid,fillcolor=gray](4,4)(5,3)
% part 4
\psframe[dimen=middle](0,3)(1,2)
\psframe[dimen=middle](1,3)(2,2)
\psframe[dimen=middle](2,3)(3,2)
\psframe[dimen=middle](3,3)(4,2)
% part 3
\psframe[dimen=middle](0,2)(1,1)
\psframe[dimen=middle](1,2)(2,1)
\psframe[dimen=middle](2,2)(3,1)
% part 3 (surlignee)
\psframe[dimen=middle](0,1)(1,0)
\psframe[dimen=middle](1,1)(2,0)
\psframe[dimen=middle,fillstyle=solid,fillcolor=gray](2,1)(3,0)
\end{pspicture}
        \caption{Ferrers diagram of the overpartition $\lambda = (\overline{5},4,3,\overline{3})$ of weight 15.}
        \label{Ferrers}
\end{figure}
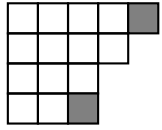

The \emph{multiplicity} of the part $j$ of an overpartition, denoted by $f_j$, is the number of occurrences of this part.
We overline the multiplicity if the part appears overlined. For example, the multiplicity of the part 4 in the overpartition $(6,6,5,4,4,\overline{4},3,\overline{1})$ is $f_4 = \overline{3}$.
The multiplicity sequence is the sequence $(f_1,f_2,\ldots )$.
For example the previous overpartition has multiplicity sequence $(\overline{1},0,1,\overline{3},1,2)$.\\

The Frobenius representation of an overpartition \cite{cl,lo2} of $n$ is  a two-rowed array
$$
\begin{pmatrix}
a_1 & a_2 & \cdots & a_N \\
b_1 & b_2 & \cdots & b_N
\end{pmatrix}
$$
where $(a_1,\ldots , a_N)$ is a partition  into distinct nonnegative
parts and  $(b_1,\ldots ,b_N)$ is an overpartition into nonnegative
parts where the {\em first} occurrence of a part can be overlined and $N +
\sum (a_i + b_i) = n$.
Following \cite{cl}, we call $p_{Q,\mathcal O}(n)$ the number of such two-rowed arrays.
We call this the Frobenius representation of an overpartition
because it is in bijection with overpartitions. This was proved in
\cite{cl} and we now recall the algorithm used for that proof.

We use the notion of a \emph{hook}. Given a positive integer $a$ and
a non-negative integer $b$, $h(a,b)$ is the hook that corresponds to
the partition $(a,1,\ldots,1)$ where there are $b$ ones. Combining a
hook $h(a,b)$ and a partition $\alpha$ is possible if and only if $a
> \alpha_1$ and $b\geq l(\alpha)$, where $l(\alpha)$ denotes the number of parts of $\alpha$. The
result of the union is $\beta = h(a,b)\cup\alpha$ with $\beta_1 =
a$, $l(\beta) = b+1$ and $\beta_i = \alpha_{i-1} +1$ for $i > 1$.

Now take a two-rowed array $\nu$ counted by $p_{Q,\mathcal O}(n)$, increase the entries on the top row by 1 and
initialize $\alpha$ and $\beta$ to the empty object, $\epsilon$.
Beginning with the rightmost column of $\nu$, we proceed to the
left, building $\alpha$ into an ordinary partition and $\beta$ into
a partition into distinct parts. At the $i^{th}$ column, if $b_i$ is
overlined, then we combine the hook $h(a_i, b_i)$ and $\alpha$.
Otherwise, we add the part $b_i$ to $\alpha'$ (the conjugate of
$\alpha$) and the part $a_i$ to $\beta$. Joining the parts of
$\alpha$ together with the parts of $\beta$ gives the overpartition
$\lambda$. An example is given below starting with $\nu =
\begin{pmatrix}7 &5 &4 &2 &0\\ 6 &\overline{4} &4 &3
&\overline{1}\end{pmatrix}$.
\begin{table}[ht]
  \centering
  \begin{tabular}{lll}
    \,\,\,$\nu$ & $\alpha$ & $\beta$ \\[0.5cm]
    $\begin{pmatrix}8 &6 &5 &3 &1\\ 6 &\overline{4} &4 &3
&\overline{1}\end{pmatrix}$ & $\epsilon$ & $\epsilon$ \\
    $\begin{pmatrix}8 &6 &5 &3\\ 6 &\overline{4} &4 &3\end{pmatrix}$ & $(1,1)$ & $\epsilon$ \\
    $\begin{pmatrix}8 &6 &5\\ 6 &\overline{4} &4\end{pmatrix}$ & $(2,2,1)$ & $(3)$ \\
    $\begin{pmatrix}8 &6\\ 6 &\overline{4}\end{pmatrix}$ & $(3,3,2,1)$ & $(5,3)$ \\
    $\begin{pmatrix}8 \\6\end{pmatrix}$
     & $(6,4,4,3,2)$ & $(5,3)$ \\
    \,\,\,\,$\epsilon$ & $(7,5,5,4,3,1)$ & $(8,5,3)$ \\
  \end{tabular}
%  \caption{}\label{}
\end{table}

We get $\lambda=
(\overline{8},7,5,5,\overline{5},4,3,\overline{3},1)$. The reverse
bijection can be easily described. See \cite{cl}.

%Given $\alpha$ and $\beta$, we set
%the Frobenius partition equal to $\epsilon$. We proceed until
%$\alpha$ and $\beta$ are empty, at each step adding a column to the
%Frobenius partition according to the following rule: If $\beta_1 \ge
%\alpha_1$ then add the column $\begin{pmatrix}\beta_1\\
%l(\alpha)\end{pmatrix}$ and decrease the parts of $\alpha$ by 1 and
%delete the largest part of $\beta$. Otherwise add the column
%$\begin{pmatrix}\alpha_1\\ l(\alpha)-1\end{pmatrix}$ and delete the
%hook $h(\alpha_1,l(\alpha)-1)$ from $\alpha$. Finally, decrease by 1
%the entries of the top row. For example, using this recipe one
%easily traces the pair $(\alpha,\beta)$ above back to the Frobenius
%partition.

We now define the successive ranks.
\begin{definition}
The \emph{successive ranks} of an overpartition can be defined from its Frobenius representation.  If an overpartition has Frobenius representation
$$
\begin{pmatrix}
a_1&a_2&\cdots&a_N\\
b_1&b_2&\cdots&b_N\\
\end{pmatrix}
$$
then its $i$th successive rank $r_i$ is $a_i - b_i$ minus the number of non-overlined parts in $\{b_{i+1},\ldots,b_N\}$.
\end{definition}
This definition is an extension of Lovejoy's definition of the rank \cite{lo2}.
For example, the successive ranks of
$
\begin{pmatrix}
7&4&2&0\\
\overline{3}&3&1&\overline{0}\\
\end{pmatrix}
$
are $(2,0,1,0)$.\\

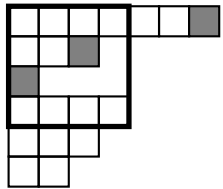
\begin{figure}[ht]
        \psset{unit=0.4}
        \begin{pspicture}(0,0)(7,6)
        % part 7
        \psframe[dimen=middle](0,6)(1,5)
        \psframe[dimen=middle](1,6)(2,5)
        \psframe[dimen=middle](2,6)(3,5)
        \psframe[dimen=middle](3,6)(4,5)
        \psframe[dimen=middle](4,6)(5,5)
        \psframe[dimen=middle](5,6)(6,5)
        \psframe[dimen=middle,fillstyle=solid,fillcolor=gray](6,6)(7,5)
        % part 3
        \psframe[dimen=middle](0,5)(1,4)
        \psframe[dimen=middle](1,5)(2,4)
        \psframe[dimen=middle,fillstyle=solid,fillcolor=gray](2,5)(3,4)
        % part 1
        \psframe[dimen=middle,fillstyle=solid,fillcolor=gray](0,4)(1,3)
        % part 4
        \psframe[dimen=middle](0,3)(1,2)
        \psframe[dimen=middle](1,3)(2,2)
        \psframe[dimen=middle](2,3)(3,2)
        \psframe[dimen=middle](3,3)(4,2)
        % part 3
        \psframe[dimen=middle](0,2)(1,1)
        \psframe[dimen=middle](1,2)(2,1)
        \psframe[dimen=middle](2,2)(3,1)
        % part 2
        \psframe[dimen=middle](0,1)(1,0)
        \psframe[dimen=middle](1,1)(2,0)
        % carre
        \psframe[dimen=middle,linewidth=2pt](0,6)(4,2)
        \end{pspicture}
        \caption{The generalized Durfee square of $\lambda=(\overline{7},4,3,\overline{3},2,\overline{1})$ has size $4$.}
        \label{exDurfeeg}
\end{figure}

We say that the \emph{generalized Durfee square} of an overpartition $\lambda$ has size $N$ if $N$ is the largest integer such that the number of overlined
parts plus the number of non-overlined parts greater or equal to $N$ is greater than
or equal to $N$ (see Figure \ref{exDurfeeg}).
\begin{proposition}
There exists a bijection between overpartitions whose Frobenius representation has $N$ columns
and whose bottom line has $j$ overlined parts and overpartitions with generalized Durfee square
of size $N$ and $N-j$ overlined parts.
\label{2.2}
\end{proposition}
\begin{proof}
An overpartition with generalized Durfee square of size $N$ can be
decomposed into a partition $\alpha$ into at most $N$ parts (the conjugate of the partition
under the generalized Durfee square) and an
overpartition $\gamma$ into $N$ parts whose non-overlined parts are
$\ge N$ (the rest). For example, the overpartition on Figure \ref{exDurfeeg} gives $\alpha=(2,2,1)$
and $\gamma=(\overline{7},4,\overline{3},\overline{1})$.
An overpartition whose Frobenius representation has $N$
columns can be decomposed into a partition $\beta$ into $N$ distinct
parts ($\beta$ is obtained by adding 1 to each part of the top
line), a partition $\delta$ into distinct parts which lie between 0
and $N-1$ and a partition $\alpha$ into at most $N$ parts.
$\delta$ and $\alpha$ are obtained from the bottom line as
follows: we first initialize $\alpha$ to the bottom line, then
if the $i$th part of the bottom line is overlined, we take off its
overline, we decrease the first $(i-1)$ parts of $\alpha$ by 1
and add a part $i-1$ to $\delta$.

Now there exists a bijection between ordered pairs $(\beta,\delta)$
and overpartitions $\gamma$. This bijection is
defined as follows: we overline all the parts of $\beta$, then for
each part $i$ in $\delta$, we add $i$ to $\beta_{i+1}$ and we remove
the overlining. We then reorder the parts, which gives us $\gamma$.
This is easily invertible and is very similar to the Algorithm Z \cite{ab2}.

The decompositions of the first paragraph and the bijection of the second paragraph 
show that there is indeed a bijection between overpartitions whose Frobenius representation has $N$ columns and overpartition with generalized Durfee square of size $N$.
If there are $j$ overlined parts in the bottom line of the Frobenius symbol, 
there are $j$ parts in $\delta$ and by the bijection, there are $N-j$ overlined parts in $\gamma$.
The proposition is thus established.
\end{proof}
\begin{example}\label{bijDF}
Let $\nu =
\begin{pmatrix}7 &5 &4 &2 &0\\ 6 &\overline{4} &4 &3
&\overline{1}\end{pmatrix}$. We thus have $N=5$ and $j=2$ in this
example. By adding 1 to each part of the top line, we get
$\beta=(8,6,5,3,1)$. From the bottom line, we get $\delta=(4,1)$ and
$\alpha=(4,3,3,2,1)$. By applying the bijection described in
the second paragraph of the above proof, we get
$\gamma=(\overline{8},7,5,\overline{5},\overline{3})$. Since
$\alpha=(4,3,3,2,1)$, the resulting overpartition is
$\lambda=(\overline{8},7,5,5,\overline{5},4,3,\overline{3},1)$. It
has $N-j=3$ overlined parts and its generalized Durfee square is of size 5.
\end{example}

This decomposition shows that the generating function for overpartitions with generalized Durfee
square of size $N$ where the exponent of $q$ counts
the weight and the exponent of $a$ counts the number of overlined parts is
$$
\frac{a^Nq^{\binom{N+1}{2}}(-1/a)_N}{(q)_N(q)_N}.
$$

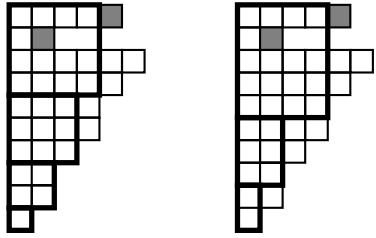
\begin{figure}[ht]
    \centering
    \psset{unit=0.3}
    \begin{pspicture}(0,0)(6,-9)
    % _5
    \psframe[dimen=middle](0,0)(1,-1)
    \psframe[dimen=middle](1,0)(2,-1)
    \psframe[dimen=middle](2,0)(3,-1)
    \psframe[dimen=middle](3,0)(4,-1)
    \psframe[dimen=middle,fillstyle=solid,fillcolor=gray](4,0)(5,-1)
    % _2
    \psframe[dimen=middle](0,-1)(1,-2)
    \psframe[dimen=middle,fillstyle=solid,fillcolor=gray](1,-1)(2,-2)
    % 6
    \psframe[dimen=middle](0,-2)(1,-3)
    \psframe[dimen=middle](1,-2)(2,-3)
    \psframe[dimen=middle](2,-2)(3,-3)
    \psframe[dimen=middle](3,-2)(4,-3)
    \psframe[dimen=middle](4,-2)(5,-3)
    \psframe[dimen=middle](5,-2)(6,-3)
    % 5
    \psframe[dimen=middle](0,-3)(1,-4)
    \psframe[dimen=middle](1,-3)(2,-4)
    \psframe[dimen=middle](2,-3)(3,-4)
    \psframe[dimen=middle](3,-3)(4,-4)
    \psframe[dimen=middle](4,-3)(5,-4)
    % 4
    \psframe[dimen=middle](0,-4)(1,-5)
    \psframe[dimen=middle](1,-4)(2,-5)
    \psframe[dimen=middle](2,-4)(3,-5)
    \psframe[dimen=middle](3,-4)(4,-5)
    % 4
    \psframe[dimen=middle](0,-5)(1,-6)
    \psframe[dimen=middle](1,-5)(2,-6)
    \psframe[dimen=middle](2,-5)(3,-6)
    \psframe[dimen=middle](3,-5)(4,-6)
    % 3
    \psframe[dimen=middle](0,-6)(1,-7)
    \psframe[dimen=middle](1,-6)(2,-7)
    \psframe[dimen=middle](2,-6)(3,-7)
    % 2
    \psframe[dimen=middle](0,-7)(2,-8)
    \psframe[dimen=middle](1,-7)(2,-8)
    % 2
    \psframe[dimen=middle](0,-8)(1,-9)
    \psframe[dimen=middle](1,-8)(2,-9)
    % Durfee squares
    \psframe[dimen=middle,linewidth=2pt](0,0)(4,-4)
    \psframe[dimen=middle,linewidth=2pt](0,-4)(3,-7)
    \psframe[dimen=middle,linewidth=2pt](0,-7)(2,-9)
    \psframe[dimen=middle,linewidth=2pt](0,-9)(1,-10)
    \end{pspicture}
    \hspace{1cm}
    \begin{pspicture}(0,0)(6,-9)
    % _5
    \psframe[dimen=middle](0,0)(1,-1)
    \psframe[dimen=middle](1,0)(2,-1)
    \psframe[dimen=middle](2,0)(3,-1)
    \psframe[dimen=middle](3,0)(4,-1)
    \psframe[dimen=middle,fillstyle=solid,fillcolor=gray](4,0)(5,-1)
    % _2
    \psframe[dimen=middle](0,-1)(1,-2)
    \psframe[dimen=middle,fillstyle=solid,fillcolor=gray](1,-1)(2,-2)
    % 6
    \psframe[dimen=middle](0,-2)(1,-3)
    \psframe[dimen=middle](1,-2)(2,-3)
    \psframe[dimen=middle](2,-2)(3,-3)
    \psframe[dimen=middle](3,-2)(4,-3)
    \psframe[dimen=middle](4,-2)(5,-3)
    \psframe[dimen=middle](5,-2)(6,-3)
    % 5
    \psframe[dimen=middle](0,-3)(1,-4)
    \psframe[dimen=middle](1,-3)(2,-4)
    \psframe[dimen=middle](2,-3)(3,-4)
    \psframe[dimen=middle](3,-3)(4,-4)
    \psframe[dimen=middle](4,-3)(5,-4)
    % 4
    \psframe[dimen=middle](0,-4)(1,-5)
    \psframe[dimen=middle](1,-4)(2,-5)
    \psframe[dimen=middle](2,-4)(3,-5)
    \psframe[dimen=middle](3,-4)(4,-5)
    % 4
    \psframe[dimen=middle](0,-5)(1,-6)
    \psframe[dimen=middle](1,-5)(2,-6)
    \psframe[dimen=middle](2,-5)(3,-6)
    \psframe[dimen=middle](3,-5)(4,-6)
    % 3
    \psframe[dimen=middle](0,-6)(1,-7)
    \psframe[dimen=middle](1,-6)(2,-7)
    \psframe[dimen=middle](2,-6)(3,-7)
    % 2
    \psframe[dimen=middle](0,-7)(2,-8)
    \psframe[dimen=middle](1,-7)(2,-8)
    % 2
    \psframe[dimen=middle](0,-8)(1,-9)
    \psframe[dimen=middle](1,-8)(2,-9)
    % Durfee rectangles
    \psframe[dimen=middle,linewidth=2pt](0,0)(4,-5)
    \psframe[dimen=middle,linewidth=2pt](0,-5)(2,-8)
    \psframe[dimen=middle,linewidth=2pt](0,-8)(1,-10)
    \end{pspicture}
    \caption{Successive Durfee squares and successive Durfee rectangles of $(6,5,\overline{5},4,4,3,2,2,\overline{2},1)$.}
    \label{sds}
\end{figure}

\begin{definition}
The \emph{successive Durfee squares} of an overpartition are its
generalized Durfee square and the successive Durfee squares of the
partition below the generalized Durfee square, if we represent the
partition as in Figure \ref{exDurfeeg}, with the overlined parts
above the non-overlined ones. We can also define similarly the
\emph{successive Durfee rectangles} by dissecting the overpartition
with $d \times (d+1)$-rectangles instead of squares (see Figure
\ref{sds}). In this case, we also impose the condition that the
partition on the right of a $d \times (d+1)$-rectangle cannot have
more than $d$ parts.
\end{definition}

These definitions imply that the generating function for overpartitions with $i-1$ successive Durfee squares followed by $k-i$ successive Durfee rectangles (the first one being a generalized Durfee square/rectangle) is
\begin{eqnarray}\label{multseriesD}
\sum_{n_1 \ge \ldots \ge n_{k-1} \ge 0} &\frac{q^{\binom{n_1+1}{2}+n_i+\ldots +n_{k-1}}(-1/a)_{n_1}a^{n_1}}{(q)_{n_1}}\nonumber\\
\times &\left(q^{n_2^2}\qbinom{n_{1}}{n_2}\right) \left(q^{n_3^2}\qbinom{n_{2}}{n_3}\right) \cdots \left(q^{n_{k-1}^2}\qbinom{n_{k-2}}{n_{k-1}}\right)
\end{eqnarray}

where
$$
\qbinom{n}{k} = \frac{(q)_n}{(q)_k (q)_{n-k}}
$$
is the generating function of partitions into at most $k$ parts less or equal to $n-k$.

\section{Paths and generating function}

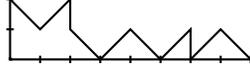
\begin{figure}[ht]
\psset{unit=0.4}
\begin{pspicture}(0,0)(7,2)
\psaxes[labels=none,ticksize=1pt](8,2)
\psline(0,2)(1,1)(2,2)(2,1)(3,0)(4,1)(5,0)(6,1)(6,0)(7,1)(8,0)
\end{pspicture}
\caption{This path has four peaks: two NES peaks (located at $(2,2)$
and $(6,1)$) and two NESE peaks (located at $(4,1)$ and $(7,1)$).
Its major index is $2+4+6+7=19$.} \label{expath}
\end{figure}

In this section we will define the paths, compute their generating function and
therefore prove Theorem \ref{thegf}. This part is an extension of papers of Andrews and Bressoud
\cite{ab,b} based on ideas of Burge \cite{bu1,bu2}. We study paths
in the first quadrant, that start on the $y$-axis, end on the $x$-axis,
and use four kinds of unitary steps:
\begin{itemize}
\item North-East $NE$: $(x,y)\rightarrow (x+1,y+1)$,
\item South-East $SE$: $(x,y)\rightarrow (x+1,y-1)$,
\item South $S$: $(x,y)\rightarrow (x,y-1)$,
\item East $E$: $(x,0)\rightarrow (x+1,0)$.
\end{itemize}
The \emph{height} corresponds to the $y$-coordinate. A South step
can only appear after a North-East step and an East step can only
appear at height 0. The paths are either empty or end with a South-East step or a South step. 
A \emph{peak} is a vertex preceded by a North-East step and
followed by a South step (in which case it will be called a
\emph{NES peak}) or by a South-East step (in which case it will be
called a \emph{NESE peak}). The \emph{major index} of a
path is the sum of the $x$-coordinates of its peaks (see Figure
\ref{expath} for an example). Let $k$ and $i$ be positive integers with $i \le k$. We say that a path satisfies the \emph{special $(k,i)$-conditions} if it starts at height $k-i$ and its height is less than $k$. When the path has no South steps, this is the definition of the paths in \cite{b}.

Let $\overline{E}_{k,i}(n,j,N)$ be the number of paths of major index $n$ with
$N$ peaks and  $j$ South steps which satisfy the special $(k,i)$-conditions.
Let $\overline{\mathcal E}_{k,i}(N)$ be the generating function for these paths, that is
$\overline{\mathcal E}_{k,i}(N)=\overline{\mathcal E}_{k,i}(N,a,q)=\sum_{n,j}\overline{E}_{k,i}(n,j,N)a^jq^n$.
Moreover, for $0 \le i <k$,  let $\overline{\Gamma}_{k,i}(N)$ be the generating function for paths obtained by deleting the first NE step of a path which is counted in $\overline{\mathcal E}_{k,i+1}(N)$ and begins with a NE step.

Then
\begin{proposition}\label{rec}
\begin{eqnarray}
\overline{\mathcal E}_{k,i}(0)&=&1;\label{rec4} \\
\overline{\mathcal E}_{k,i}(N)&=&q^N\overline{\mathcal E}_{k,i+1}(N)+q^N\overline{\Gamma}_{k,i-1}(N);\ \ \ \ i<k \label{rec1}\\
\overline{\Gamma}_{k,i}(N)&=&q^N\overline{\Gamma}_{k,i-1}(N)+(a+q^{N-1})\overline{\mathcal E}_{k,i+1}(N-1);\ \ \ i>0 \label{rec2}\\
\overline{\mathcal E}_{k,k}(N)&=& \frac{q^N}{1-q^N} \overline{\Gamma}_{k,k-1}(N); \label{rec3}\\
 \overline{\Gamma}_{k,0}(N)&=&0.\label{rec5}
\end{eqnarray}
\end{proposition}
\begin{proof}
If $N=0$ the only path counted in $\overline{\mathcal E}_{k,i}(0)$ is the path with only South-East steps, which starts at $(0,k-i)$ and ends at $(k-i,0)$ (if $i=k$, it is just the empty path, starting and ending at $(0,0)$).
This path has no peaks and its major index is thus 0. This proves \eqref{rec4}.
Now if the path has at least one peak, then we take off its first step.
If $i<k$, then a path counted in $\overline{\mathcal E}_{k,i}(N)$ starts with a North-East
(defined by $q^N\overline{\Gamma}_{k,i-1}(N)$) or a South-East step ($q^N\overline{\mathcal E}_{k,i+1}(N)$). This gives \eqref{rec1}.
If $i>0$, $\overline{\Gamma}_{k,i}(N)$ is the generating function for paths counted by $\overline{\mathcal E}_{k,i+1}(N)$ that start with a North-East step
where the first step was
deleted. These paths can start with a North-East step ($q^N\overline{\Gamma}_{k,i-1}(N)$), a South
step ($a\overline{\mathcal E}_{k,i+1}(N-1)$) or a South-East step ($q^{N-1}\overline{\mathcal E}_{k,i+1}(N-1)$) and we get \eqref{rec2}.
If $i=k$ then a path counted by $\overline{\mathcal E}_{k,k}(N)$ starts with a North-East ($q^N\overline{\Gamma}_{k,k-1}(N)$) or an East step ($q^N\overline{\mathcal E}_{k,k}(N)$).
This gives \eqref{rec3}.
The height of the paths is less than $k$, therefore no path which starts at height $k-1$ can start
with  a North-East step and $\overline{\Gamma}_{k,0}(N)=0$ as in \eqref{rec5}.
\end{proof}

These recurrences uniquely define the series $\overline{\mathcal
E}_{k,i}(N)$ and $\overline{\Gamma}_{k,i}(N)$. We get that:
\begin{theorem}
\begin{eqnarray*}
\overline{\mathcal E}_{k,i}(N)&=&a^Nq^{\binom{N+1}{2}}(-1/a)_N\sum_{n=-N}^{N} (-1)^n \frac{q^{kn^2+n(k-i)-{\binom{n}{2}}}}{(q)_{N-n}(q)_{N+n}}\\
\overline{\Gamma}_{k,i}(N)&=&a^Nq^{\binom{N}{2}}(-1/a)_N\sum_{n=-N}^{N-1} (-1)^n \frac{q^{kn^2+n(k-i)-{\binom{n+1}{2}}}}{(q)_{N-n-1}(q)_{N+n}}
\end{eqnarray*}
\end{theorem}

\begin{proof}
Let
\begin{equation*}
\left\{
\begin{aligned}
\overline{\mathcal E'}_{k,i}(N)&=a^Nq^{\binom{N+1}{2}}(-1/a)_N\sum_{n=-N}^{N} (-1)^n \frac{q^{kn^2+n(k-i)-{\binom{n}{2}}}}{(q)_{N-n}(q)_{N+n}} \\
\overline{\Gamma'}_{k,i}(N)&=a^Nq^{\binom{N}{2}}(-1/a)_N\sum_{n=-N}^{N-1} (-1)^n \frac{q^{kn^2+n(k-i)-{\binom{n+1}{2}}}}{(q)_{N-n-1}(q)_{N+n}}
\end{aligned}
\right.
\end{equation*}
Note that $\overline{\mathcal E'}_{k,i}(0) = 1$.

We first prove that $\overline{\mathcal E'}_{k,i}(N)$ and $\overline{\Gamma'}_{k,i}(N)$ satisfy $\overline{\mathcal E'}_{k,i}(N)=q^N\overline{\mathcal E'}_{k,i+1}(N) +
q^N \overline{\Gamma'}_{k,i-1}(N)$ for $1\le i\le k$ :
\begin{align*}
&q^N\overline{\mathcal E'}_{k,i+1}(N) + q^N \overline{\Gamma'}_{k,i-1}(N)\\
=& a^Nq^{\binom{N+1}{2}}(-1/a)_N\sum_{n=-N}^{N} (-1)^n \frac{q^{kn^2+n(k-i-1)-{\binom{n}{2}}}}{(q)_{N-n}(q)_{N+n}} q^N\\
& + a^Nq^{\binom{N}{2}}(-1/a)_N\sum_{n=-N}^{N-1} (-1)^n \frac{q^{kn^2+n(k-i+1)-{\binom{n+1}{2}}}}{(q)_{N-n-1}(q)_{N+n}} q^{N}\\
=& a^Nq^{\binom{N+1}{2}}(-1/a)_N\sum_{n=-N}^{N-1} (-1)^n \frac{q^{kn^2+n(k-i)-{\binom{n}{2}}}}{(q)_{N-n}(q)_{N+n}} (q^{N-n}+(1-q^{N-n}))\\
& + a^Nq^{\binom{N+1}{2}}(-1/a)_N (-1)^N \frac{q^{kN^2+N(k-i-1)-{\binom{N}{2}}}}{(q)_{0}(q)_{2N}} q^N\\
=& a^Nq^{\binom{N+1}{2}}(-1/a)_N\sum_{n=-N}^{N} (-1)^n \frac{q^{kn^2+n(k-i)-{\binom{n}{2}}}}{(q)_{N-n}(q)_{N+n}}\\
=& \overline{\mathcal E'}_{k,i}(N)
\end{align*}

We then prove, quite similarly, that they satisfy \eqref{rec2}:
\begin{align*}
&(a+q^{N-1}) \overline{\mathcal E'}_{k,i+1}(N-1) + q^N \overline{\Gamma'}_{k,i-1}(N)\\
= &a^{N-1}q^{\binom{N}{2}}(-1/a)_{N-1}\sum_{n=-N+1}^{N-1} (-1)^n \frac{q^{kn^2+n(k-i-1)-{\binom{n}{2}}}}{(q)_{N-1-n}(q)_{N-1+n}}(a+q^{N-1})\\
& + a^Nq^{\binom{N}{2}}(-1/a)_N\sum_{n=-N}^{N-1} (-1)^n \frac{q^{kn^2+n(k-i+1)-{\binom{n+1}{2}}}}{(q)_{N-n-1}(q)_{N+n}} q^N\\
= & a^Nq^{\binom{N}{2}}(-1/a)_N\sum_{n=-N+1}^{N-1} (-1)^n \frac{q^{kn^2+n(k-i)-{\binom{n+1}{2}}}}{(q)_{N-n-1}(q)_{N+n}} ((1-q^{N+n})+q^{N+n})\\
& + a^Nq^{\binom{N}{2}}(-1/a)_N (-1)^{-N} \frac{q^{kN^2-N(k-i+1)-{\binom{-N+1}{2}}}}{(q)_{2N-1}(q)_{0}} q^N\\
= & a^Nq^{\binom{N}{2}}(-1/a)_N\sum_{n=-N}^{N-1} (-1)^n \frac{q^{kn^2+n(k-i)-{\binom{n+1}{2}}}}{(q)_{N-n-1}(q)_{N+n}}\\
= & \overline{\Gamma'}_{k,i}(N)
\end{align*}

For \eqref{rec3}, we prove that $\overline{\mathcal E'}_{k,k+1}(N) =
\overline{\mathcal E'}_{k,k}(N)$:
\begin{align*}
\overline{\mathcal E'}_{k,k+1}(N) = &a^Nq^{\binom{N+1}{2}}(-1/a)_N\sum_{n=-N}^{N} (-1)^n \frac{q^{kn^2-n-{\binom{n}{2}}}}{(q)_{N-n}(q)_{N+n}}\\
= & a^Nq^{\binom{N+1}{2}}(-1/a)_N\sum_{p=-N}^{N} (-1)^{-p} \frac{q^{k(-p)^2+p-{\binom{-p}{2}}}}{(q)_{N+p}(q)_{N-p}} \intertext{ where $p=-n$}\\
= &a^Nq^{\binom{N+1}{2}}(-1/a)_N\sum_{p=-N}^{N} (-1)^p \frac{q^{kp^2+p-\frac{(-p)(-p-1)}{2}}}{(q)_{N+p}(q)_{N-p}}\\
= &a^Nq^{\binom{N+1}{2}}(-1/a)_N\sum_{p=-N}^{N} (-1)^p \frac{q^{kp^2+p-\frac{(p)(p+1)}{2}}}{(q)_{N+p}(q)_{N-p}}\\
= &a^Nq^{\binom{N+1}{2}}(-1/a)_N\sum_{p=-N}^{N} (-1)^p \frac{q^{kp^2+p-\binom{p+1}{2}}}{(q)_{N+p}(q)_{N-p}}\\
= &a^Nq^{\binom{N+1}{2}}(-1/a)_N\sum_{p=-N}^{N} (-1)^p \frac{q^{kp^2-\binom{p}{2}}}{(q)_{N+p}(q)_{N-p}}\\
= & \overline{\mathcal E'}_{k,k}(N)
\end{align*}
Hence we have, using the fact that $\overline{\mathcal E'}_{k,i}(N)$ satisfies \eqref{rec1} for $i=k$:
$$
\overline{\mathcal E'}_{k,k}(N) = q^N \overline{\mathcal E'}_{k,k}(N) + q^N \overline{\Gamma'}_{k,k-1}(N).
$$
\par Finally, we have:
\begin{align*}
&\overline{\Gamma'}_{k,0}(N) = a^Nq^{\binom{N}{2}}(-1/a)_N\sum_{n=-N}^{N-1} (-1)^n \frac{q^{kn^2+nk-{\binom{n+1}{2}}}}{(q)_{N-n-1}(q)_{N+n}}\\
= & a^Nq^{\binom{N}{2}}(-1/a)_N \left( \sum_{n=0}^{N-1} (-1)^n \frac{q^{kn^2+nk-{\binom{n+1}{2}}}}{(q)_{N-n-1}(q)_{N+n}}
+ \sum_{n=-N}^{-1} (-1)^n \frac{q^{kn^2+nk-{\binom{n+1}{2}}}}{(q)_{N-n-1}(q)_{N+n}} \right)\\
= & a^Nq^{\binom{N}{2}}(-1/a)_N \left( \sum_{n=0}^{N-1} (-1)^n \frac{q^{kn^2+nk-{\binom{n+1}{2}}}}{(q)_{N-n-1}(q)_{N+n}}\right.\\
&+ \left.\sum_{p=0}^{N-1} (-1)^{-1-p} \frac{q^{k(-1-p)^2+(-1-p)k-{\binom{-p}{2}}}}{(q)_{N+p}(q)_{N-1-p}} \right) \intertext{where $p=-1-n$}\\
= & a^Nq^{\binom{N}{2}}(-1/a)_N \left( \sum_{n=0}^{N-1} (-1)^n \frac{q^{kn^2+nk-{\binom{n+1}{2}}}}{(q)_{N-n-1}(q)_{N+n}}
- \sum_{p=0}^{N-1} (-1)^p \frac{q^{kp^2 +kp-{\binom{p+1}{2}}}}{(q)_{N+p}(q)_{N-1-p}} \right)\\
= & 0
\end{align*}

Since $\overline{\mathcal E'}_{k,i}(N)$ and $\overline{\Gamma'}_{k,i}(N)$ satisfy the recurrences
of Proposition \ref{rec}, we thus have $\overline{\mathcal E}_{k,i}(N)=\overline{\mathcal E'}_{k,i}(N)$ and $\overline{\Gamma}_{k,i}(N)=\overline{\Gamma'}_{k,i}(N)$.
\end{proof}

We just need the following proposition to prove Theorem \ref{thegf}.
\begin{proposition} \label{superDurfee}
For any $n \in \mathbb{Z}$
$$
\sum_{N\ge |n|}
\frac{(-aq)_n(-q^n/a)_{N-n}q^{\binom{N+1}{2}-\binom{n+1}{2}}a^{N-n}}{(q)_{N+n}(q)_{N-n}}=
\frac{(-aq)_\infty}{(q)_\infty}.
$$
\end{proposition}
\begin{proof}
We only prove the case $n\ge 0$. The case $n<0$ is identical as $(a)_{-n}=(-1/a)^nq^{n+1\choose 2}/(aq)_n$ and therefore
$$
\frac{(-aq)_n(-q^n/a)_{N-n}q^{\binom{N+1}{2}-\binom{n+1}{2}}a^{N-n}}{(q)_{N+n}(q)_{N-n}}=
\frac{(-aq)_{-n}(-q^{-  n}/a)_{N+n}q^{\binom{N+1}{2}-\binom{-n+1}{2}}a^{N+n}}{(q)_{N+n}(q)_{N-n}}.
$$
We present an analytical proof and a combinatorial one.

{\em Analytical proof: }
\begin{align*}
&\sum_{N=n}^{\infty} \frac{(-aq)_n(-q^n/a)_{N-n}q^{\binom{N+1}{2}-\binom{n+1}{2}}a^{N-n}}{(q)_{N+n}(q)_{N-n}} \\
= &\sum_{N =0}^{\infty} \frac{(-aq)_n(-q^n/a)_N q^{\binom{N+n+1}{2}-\binom{n+1}{2}}a^N}{(q)_{N}(q)_{N+2n}}\\
= &\frac{(-aq)_n}{(q)_{2n}} \sum_{N=0}^{\infty}\frac{q^{Nn+\binom{N+1}{2}}(-q^n/a)_N a^N}{(q)_N(q^{2n+1})_N}\\
\intertext{We now apply the $q$-Gauss summation (Corollary 2.4 of \cite{ag} with $n\to N$, $a \to -q^n/a$, $b \to -\infty$ and $c \to q^{2n+1}$)}
= &\frac{(-aq)_n}{(q)_{2n}} \frac{(-aq^{n+1})_\infty}{(q^{2n+1})_\infty}\\
= &\frac{(-aq)_\infty}{(q)_\infty}
\end{align*}

{\em Combinatorial proof:} Let $n$ be a fixed non-negative
integer and let $N$ be the greatest integer such that the
sum of the number of overlined parts greater than $n$ and of the number of
non-overlined parts greater than or equal to $N+n$ is greater than or equal to $N-n$.
By definition, $N$,
if it exists,  is unique. We check that $N=n$ satisfies the condition: the
sum of the number of overlined parts $>n$ and of the number of
non-overlined parts $\geq 2n$ is indeed $\geq 0$. Therefore, $N$
exists and it is unique. We call $N$ the size of the generalized $n$-Durfee square.
Note that if $n=0$, $N$ is the size of the generalized Durfee square and this bijection is
the same as the one presented in the proof of Proposition \ref{2.2}.
For example, if $n=2$ then the overpartition $(8,8,\bar{6},5,5,3,3,\bar{3},\bar{1})$
gives $N=6$.

We now show that the generating function for the overpartitions with generalized
$n$-Durfee square $N$ is:
$$
\frac{(-aq)_n(-q^n/a)_{N-n}q^{\binom{N+1}{2}-\binom{n+1}{2}}a^{N-n}}{(q)_{N+n}(q)_{N-n}}.
$$

The factor $(-aq)_n$ corresponds to the overlined parts $\leq n$
and the factor $\frac{1}{(q)_{N+n}}$ corresponds to the
non-overlined parts $\leq N+n$.

The remaining factors correspond to an overpartition into $N-n$
parts whose overlined parts are $>n$ and whose non-overlined parts
are $\geq N+n$. To prove this, let us show that there exists a
bijection between such overpartitions $\theta$ and triples
$(\varepsilon,\zeta,\eta)$ where $\varepsilon$ is 
the overpartition $(\overline{N},\overline{N-1},\ldots ,\overline{n+1})$
which has $N-n$ parts
(corresponding to the factor $a^{N-n}q^{\binom{N+1}{2}-\binom{n+1}{2}}$), $\zeta$ is a
partition into distinct parts which lie between $n$ and $N-1$
(corresponding to the factor $(-q^n/a)_{N-n}$) and $\eta$ is a partition into
$N-n$ non-negative parts (corresponding to the factor $\frac{1}{(q)_{N-n}}$). This bijection (similar
to Algorithm Z \cite{ab2}) is defined as follows: first for all $i$ we
set $\theta_i=\varepsilon_i+\eta_i$. Then for each part $n+i$ in $\zeta$, we add $n+i$ to $\theta_{i+1}$
and we remove the overlining of that part. This implies that the non-overlined parts
are $\ge N+n$. Finally, we reorder the parts. It is easy to see that this is a bijection as there
only a unique ordering of the parts of $\theta$ which allows,
if $\theta_{i+1}$ is not overlined,  to take off $n+i$ from it and overline it
and get a partition into distinct overlined parts.

%\begin{figure}[ht]
%    \centering
 %   \psset{unit=0.4}
 %   \begin{pspicture}(0,0)(7,-4)
  %  %\pspolygon(0,0)(7,0)(7,-1)(6,-1)(6,-2)(5,-2)(5,-3)(4,-3)(4,-4)(0,-4)
    % 7
  %  \psframe[dimen=middle](0,0)(1,-1)
    %\psframe[dimen=middle](1,0)(2,-1)
    %\psframe[dimen=middle](2,0)(3,-1)
    %\psframe[dimen=middle](3,0)(4,-1)
    %\psframe[dimen=middle](4,0)(5,-1)
    %\psframe[dimen=middle](5,0)(6,-1)
    %\psframe[dimen=middle](6,0)(7,-1)
    % 6
    %\psframe[dimen=middle](0,-1)(1,-2)
    %\psframe[dimen=middle](1,-1)(2,-2)
    %\psframe[dimen=middle](2,-1)(3,-2)
    %\psframe[dimen=middle](3,-1)(4,-2)
    %\psframe[dimen=middle](4,-1)(5,-2)
    %\psframe[dimen=middle](5,-1)(6,-2)
    % 5
    %\psframe[dimen=middle](0,-2)(1,-3)
    %\psframe[dimen=middle](1,-2)(2,-3)
    %\psframe[dimen=middle](2,-2)(3,-3)
    %\psframe[dimen=middle](3,-2)(4,-3)
    %\psframe[dimen=middle](4,-2)(5,-3)
    % 4
    %\psframe[dimen=middle](0,-3)(1,-4)
    %\psframe[dimen=middle](1,-3)(2,-4)
    %\psframe[dimen=middle](2,-3)(3,-4)
    %\psframe[dimen=middle](3,-3)(4,-4)
    %\end{pspicture}
    %\caption{The partition $\varepsilon$ for $N=7$ and $n=3$.}
    %\label{trapezoid}
%\end{figure}

For example for $N=7$ and $n=3$, if $\varepsilon=(\bar{7},\bar{6},\bar{5},\bar{4})$, $\eta=(3,3,1,0)$
and $\zeta=(6,3)$, then $\theta=(13,\bar{9},\bar{6},10)=(13,10,\bar{9},\bar{6})$.

%Let us show that the decomposition used in this proof is unique. 

\end{proof}

Summing on $N$ we get
\begin{eqnarray*}
\sum_{N\ge 0}\overline{\mathcal E}_{k,i}(N) &=& \sum_{N\ge 0} a^Nq^{\binom{N+1}{2}}(-1/a)_N\sum_{n=-N}^{N} (-1)^n \frac{q^{kn^2+n(k-i)-{\binom{n}{2}}}}{(q)_{N-n}(q)_{N+n}}\\
&=& \sum_{N\ge 0} \sum_{n=-N}^{N} a^Nq^{\binom{N+1}{2}}(-1/a)_N (-1)^n \frac{q^{kn^2+n(k-i)-{\binom{n}{2}}}}{(q)_{N-n}(q)_{N+n}}\\
&=& \sum_{n=-\infty}^{\infty} \sum_{N \ge |n|} a^Nq^{\binom{N+1}{2}}(-1/a)_N (-1)^n \frac{q^{kn^2+n(k-i)-{\binom{n}{2}}}}{(q)_{N-n}(q)_{N+n}}\\
&=& \sum_{n=-\infty}^{\infty} (-1)^n a^n q^{kn^2+(k-i+1)n}
\frac{(-1/a)_n}{(-aq)_n} \cdot\\
&&\quad \quad \quad \quad \sum_{N \ge |n|}
\frac{(-aq)_n(\frac{-q^n}{a})_{N-n}
q^{\binom{N+1}{2}-\binom{n+1}{2}} a^{N-n}}{(q)_{N-n}(q)_{N+n}}\\
&=& \sum_{n=-\infty}^{\infty} (-1)^n a^n q^{kn^2+(k-i+1)n}
\frac{(-1/a)_n}{(-aq)_n} \frac{(-aq)_\infty}{(q)_\infty}\\ &&\text{
by
Proposition \ref{superDurfee}}\\
&=& \frac{(-aq)_\infty}{(q)_\infty}\sum_{n=-\infty}^{\infty}(-1)^n
a^n q^{kn^2+(k-i+1)n}\frac{(-1/a)_n}{(-aq)_n}.
\end{eqnarray*}
This is equation (\ref{gf}) of Theorem \ref{thegf}.

\section{Paths and successive ranks}

In this section we prove the case $\overline{C}_{k,i}(n,j)=\overline{E}_{k,i}(n,j)$ of Theorem \ref{main}.
%This section is a generalization of Bressoud's correspondence for partitions presented in \cite{b} and based on ideas of Burge \cite{bu1,bu2}.
In fact, we prove a refinement of this case:
\begin{proposition}
There exists a one-to-one correspondence between the paths of major index $n$ with $j$ south steps
counted by $\overline{E}_{k,i}(n,j)$
and  the overpartitions of $n$ with $j$ non-overlined parts in the bottom line of their Frobenius representation and
whose successive ranks lie in $[-i+2,2k-i-1]$
counted by $\overline{C}_{k,i}(n,j)$. This correspondence is such that the paths
have $N$ peaks if and only if the Frobenius representation of the overpartition has $N$ columns.
\label{ranks}
\end{proposition}

This proposition can be proved with a recursive argument.
We can show that $\overline{\mathcal C}_{k,i}(N)$, the generating function for overpartitions
whose Frobenius representation has $N$ columns and whose
successive ranks lie in $[-i+2,2k-i-1]$, follows the same recurrences as
$\overline{\mathcal E}_{k,i}(N)$, as done by Burge \cite{bu1,bu2} for the case of partitions.
Details are given in \cite{m}.
We propose here a direct mapping that is a generalization of a bijection of Bressoud  \cite{b}.\\

Given a lattice path which starts at $(0,a)$ and a peak $(x,y)$ with $u$ South steps to its left, we map this peak to the pair $(s,t)$ where
\begin{align*}
s &= (x+a-y+u)/2 \\
t &= (x-a+y-2-u)/2
\end{align*}
if there are an even number of East steps to the left of the peak (we then say that the peak is of \emph{type 0}), and
\begin{align*}
s &= (x+a+y-1+u)/2 \\
t &= (x-a-y-1-u)/2
\end{align*}
if there are an odd number of East steps to the left of the peak (we then say that the peak is of \emph{type 1}).
Moreover, we overline $t$ if the peak is a NESE peak.
In both cases, $s$ and $t$ are integers and we have $s+t+1=x$.
In the case of partitions treated in \cite{b}, $u$ is always 0.

Let $N$ be the number of peaks in the path and $j$ the number of South steps of the paths.
Let $(x_i,y_i)$ be the coordinates of the $i^{th}$ peak from the right and $(s_i,t_i)$ be the corresponding pair.
\begin{proposition}
The sequence $(s_1,s_2,\ldots,s_N)$ is a partition into distinct nonnegative parts and
the sequence $(t_1,t_2,\ldots,t_N)$ is an overpartition into nonnegative parts with $j$ non-overlined parts.
\end{proposition}
\begin{proof}
We need to prove the following results:
\begin{itemize}
\item \emph{$s_N\ge 0$}.
If the last peak from the right $(x_N,y_N)$ is of type 0,
then $s_N=(x_N+a-y_N+u_N)/2$ and it is sufficient to prove that $x_N-y_N
\ge -a$ since $u_N= 0$. It is obvious that any vertex has a greater
(or equal) value of $x-y$ than the previous vertex in the path.
Since the path begins at $(0,a)$, we have $x-y=-a$ at the beginning
of the path and thus we have $x-y \ge -a$ for all vertices. Now if the peak is of type 1, then
$x_N\ge 2$ and $s_N=(x_N+a+y_N-1+u_N)/2\ge 0$.
% For all $i$, we have $x_i + y_i > x_{i+1}
%+ y_{i+1}$ and $u_i \ge u_{i+1}$, so $x_i+y_i+u_i
%> x_{i+1}+y_{i+1}+u_{i+1}$. Let $p$ be the leftmost peak of type 1
%in the path. We have $x_p \ge a+1+y_p+u_p$ (see Figure
%\ref{firstpeaktype1}), $y_p \ge 1$, $u_p \ge 0$ and thus $s_p \ge
%0$. Hence all the $s_i$ are nonnegative for peaks of type 1.
% remarque: les figures firstpeaktype1 et firstpeak ont été inversées
%\begin{figure}[ht]
%\centering \psset{unit=0.6}
%\begin{pspicture}(0,-1)(10,4)
%\psaxes[labels=none,ticks=none](10,4)
%\psline(0,4)(3,1)(4,2)(6,0)(7,0)(9,2)(10,1)
%\psline[arrows=<->](0,4)(6,4) \psline[linestyle=dashed](6,0)(6,4)
%\psline[arrows=<->](7,2)(9,2) \psline[linestyle=dashed](7,0)(7,2)
%\psline[linestyle=dashed](9,2)(9,0) \rput[t](3,3.75){$\ge a+u_p$}
%\rput[b](8,2.25){$\ge y_p$} \rput[t](9,-0.25){$x_p$}
%\end{pspicture}
%\caption{If $p$ is the leftmost peak of type 1, $x_p \ge
%a+1+y_p+u_p$.} \label{firstpeaktype1}
%\end{figure}
\item \emph{$t_N\ge 0$}. If the last peak from the right $(x_N,y_N)$ is of type 0, then
$t_N=(x_N-a+y_N-2-u_N)/2$.
We have
$u_N=0$ and $x_N>0$. If $y_N\ge a$, then $t_N\ge 0$. Otherwise
it is easy to see that $x_N\ge 2+a-y_N$
(see Figure \ref{firstpeak}) and thus $t_N \ge 0$.
Now if this peak is of type 1, then $t_N=(x_N-a-y_N-1-u_N)/2$. We have $u_N=0$ and
there is at least one East step to the left of the peak and therefore at least
$a$ South-East steps before the East step and at least $y_N$ North-East steps
after the East step. Hence $x_N-y_N-a>0$ and $t_N\ge 0$.
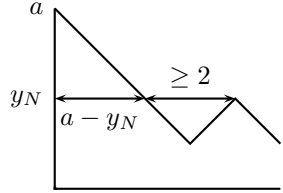
\begin{figure}[ht]
\centering \psset{unit=0.6}
\begin{pspicture}(0,0)(5,4)
\psaxes[labels=none,ticks=none](5,4) \psline(0,4)(3,1)(4,2)(5,1)
\psline[arrows=<->](0,2)(2,2) \psline[arrows=<->](2,2)(4,2)
\rput[t](1,1.75){$a-y_N$} \rput[b](3,2.25){$\ge 2$}
\rput[r](-0.25,4){$a$} \rput[r](-0.25,2){$y_N$}
\end{pspicture}
\caption{We show that $x_N \ge 2+a-y_N$ in the case $y_N \le a$.}
\label{firstpeak}
\end{figure}
\item \emph{The sequence $s$ is a partition into distinct parts}.
We need to prove that for all $i$, $s_i - s_{i+1} > 0$.
If the $i^{th}$ peak from the right and the $(i+1)^{st}$ peak are both of type 0, it is clearly true since $x_i-y_i > x_{i+1}-y_{i+1}$ (two peaks cannot have the same value of $x-y$), and $u_i \ge u_{i+1}$ (remember that $u_i$ is the number of South steps to the left of the $i^{th}$ peak).
If the $i^{th}$ peak from the right and the $(i+1)^{st}$ peak  are both of type 1,
then $s_i-s_{i+1}=(x_i+y_i+u_i-x_{i+1}-y_{i+1}-u_{i+1})/2$.
It is easy to see that that $x_i+y_i> x_{i+1}+y_{i+1}$ and $u_i\ge u_{i+1}$.
If the $i^{th}$ peak is of type 0 and the $(i+1)^{st}$ is of type 1,
we have $s_i-s_{i+1} = \frac{1}{2} (x_i-x_{i+1}-y_i-y_{i+1}+u_i-u_{i+1}+1)$.
Since $u_i \ge u_{i+1}$, it is sufficient to prove that $x_i-x_{i+1}-y_i-y_{i+1} \ge 0$.
This comes from the fact that there is an East step between the two peaks (see Figure \ref{s01}).
If the $i^{th}$ peak is of type 1 and the $(i+1)^{st}$ is of type 0, the proof is similar.
\begin{figure}[ht]
\centering \psset{unit=0.6}
\begin{pspicture}(0,-1)(8,4)
\psline(0,0)(10,0) \psline(0,4)(4,0)(5,0)(8,3)
\psline[arrows=<->](0,4)(4,4) \psline[arrows=<->](5,3)(8,3)
\psline[linestyle=dashed](0,4)(0,0)
\psline[linestyle=dashed](4,4)(4,0)
\psline[linestyle=dashed](5,3)(5,0)
\psline[linestyle=dashed](8,3)(8,0) \rput[t](2,3.75){$\ge y_{i+1}$}
\rput[t](6.5,3.75){$\ge y_i$} \rput[t](0,-0.25){$x_{i+1}$}
\rput[t](8,-0.25){$x_i$}
\end{pspicture}
\caption{If the $i^{th}$ peak of is type 0 and the $(i+1)^{st}$  peak is of type 1,
we have $x_i-x_{i+1} \ge y_i + y_{i+1}$.} \label{s01}
\end{figure}
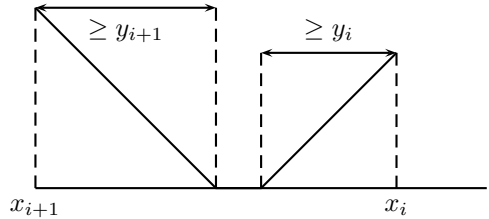
\item \emph{The sequence $t$ is an overpartition (where the first
occurrence of a part can be overlined)}.
We need to prove that for all $i$, $t_i - t_{i+1} > 0$ if $t_{i+1}$
is overlined and $t_i - t_{i+1} \ge 0$ otherwise. The fact that we
always have $t_i - t_{i+1} \ge 0$ is proved in the same way as with
the $s_i$. If $t_{i+1}$ is overlined, then it corresponds to a NESE
peak, so we have $x_i - x_{i+1} \ge 2$ and $u_{i+1} = u_i$. By
considering the expression of $t_i - t_{i+1}$ in the four cases (the peaks $i$
and $i+1$ are both of type 0,  both of type 1,  of type 0 and of type 1, or of type 1 and of type
0), the result is easily shown.
\end{itemize}
\end{proof}

Therefore $\begin{pmatrix}s_1&s_2&\cdots&s_N\\ t_1&t_2&\cdots&t_N \end{pmatrix}$ is the Frobenius representation of an overpartition whose weight is
$$
\sum_{i=1}^N (s_i+t_i+1) = \sum_{i=1}^N x_i
$$
i.e.\ the major index of the corresponding path.

\begin{figure}[ht]
\psset{unit=0.4}
\begin{pspicture}(0,0)(22,6)
\psaxes[labels=none,ticksize=1pt](0,0)(25,4)
\psline(0,2)(2,0)(6,4)(8,2)(9,3)(9,2)(10,1)(12,3)(15,0)(16,0)(18,2)(19,1)(22,4)(22,3)(25,0)
\psdots[dotstyle=x](1,1)
\psdots[dotstyle=x](3,1)
\psdots[dotstyle=x](4,2)
\psdots[dotstyle=x](5,3)
\psdots[dotstyle=x](7,3)
\psdots[dotstyle=x](11,2)
\psdots[dotstyle=x](13,2)
\psdots[dotstyle=x](14,1)
\psdots[dotstyle=x](17,1)
\psdots[dotstyle=x](20,2)
\psdots[dotstyle=x](21,3)
\psdots[dotstyle=x](23,2)
\psdots[dotstyle=x](24,1)
\rput[b](6,5){\footnotesize$(6,4,0)$}
\rput[b](9,4){\footnotesize$(9,3,0)$}
\rput[b](12,4){\footnotesize$(12,3,1)$}
\rput[b](18,3){\footnotesize$(18,2,1)$}
\rput[b](22,5){\footnotesize$(22,4,1)$}
\end{pspicture}
\caption{Illustration of the correspondence between paths and successive ranks. The values of $x$, $y$ and $u$ are given for each peak.}
\label{expathranks}
\end{figure}
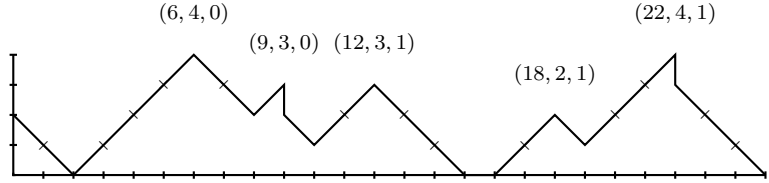

As an example, the path in Figure \ref{expathranks} corresponds to the overpartition
$$
\begin{pmatrix}
14 &11 &6 &4 &2\\
7 &\overline{6} &\overline{5} &4 &\overline{3}
\end{pmatrix}.
$$

The peaks all have height at least one, thus for a peak $(x,y)$
which is preceded by an even number of East steps, we have:
\begin{align*}
&1 \le y=a+1+t-s+u\le k-1 \\
\Leftrightarrow~&a-k+2\le s-t-u \le a \\
\Leftrightarrow~&\text{the corresponding successive rank is
$\ge a-k+2$ and $\le a$}
\end{align*}
and if the peak is preceded by an odd number of East steps, we have:
\begin{align*}
&1 \le y=s-t-u-a\le k-1 \\
\Leftrightarrow~&a+1\le s-t-u \le k+a-1 \\
\Leftrightarrow~&\text{the corresponding successive rank is $\ge a+1$
and $\le k+a-1$}.
\end{align*}
Thus, given a Frobenius representation of an overpartition and a nonnegative integer $a$, there is a unique corresponding path which starts at $(0,a)$.

In our paths, $a=k-i$, therefore in the first case the successive rank $r \in [-i+2,k-i]$ and in the second case $r \in [k-i+1,2k-i-1]$.

The map is easily reversible. This proves Proposition \ref{ranks}.

\section{Paths and multiplicities}

In this section we prove the case $\overline{B}_{k,i}(n,j)=\overline{E}_{k,i}(n,j)$ of Theorem \ref{main}.
We even prove a refinement:
\begin{proposition}
There exists a one-to-one correspondence between the paths counted by $\overline{E}_{k,i}(n,j)$
and  the overpartitions counted by $\overline{B}_{k,i}(n,j)$. This correspondence is such that the paths
have $N$ peaks if and only if the length of the multiplicity sequence of the overpartition is $N$ (see Section \ref{combproofmult} for the definition).
\label{sec5}
\end{proposition}
We will first give a generating function proof of that proposition (without the refinement).
Then we will give the sketch of a combinatorial proof which is a generalization of Burge's correspondence for partitions presented in \cite{bu1}. 

Recall that $\overline{B}_{k,i}(n,j)$ is the number of overpartitions $\lambda$ of $n$ with $j$ overlined parts such that
for all $\ell$,
\begin{equation*}
\begin{array}{ll}
&f_1 < i\\
&\lambda_{\ell}-\lambda_{\ell+k-1} \ge
\begin{cases}
1 & \text{if $\lambda_{\ell+k-1}$ is overlined}\\
2 & \text{otherwise}
\end{cases}
\end{array}
\end{equation*}
or equivalently,
\begin{equation*}
\begin{array}{ll}&f_1 < i\\
&\forall \ell \text{, } f_{\ell} + f_{\ell+1} <
\begin{cases}
k+1& \text{if a part $\ell$ is overlined}\\
k& \text{otherwise}
\end{cases}\\
\end{array}
\end{equation*}
We will abbreviate this last condition with the notation
$\forall \ell \text{, } f_{\ell} + f_{\ell+1}<k\overline{+1}$.

\subsection{A generating function proof}
Let $\overline{\mathcal B}_{k,i}(a,q)=\sum_{n,j\ge 0} \overline{B}_{k,i}(n,j) a^j q^n$.
We prove that
\begin{equation}
\overline{\mathcal B}_{k,i}(a,q) = \overline{\mathcal E}_{k,i}(a,q)
\label{mk}
\end{equation}

We will generalize Lovejoy's proof of Theorem 1.1 of \cite{lo}. Let
\begin{eqnarray}
J_{k,i}(a,x,q)& = &H_{k,i}(a,xq,q) - axq H_{k,i-1}(a,xq,q)\label{defJ}\\
H_{k,i}(a,x,q)&=&\sum_{n=0}^{\infty} \frac{x^{kn}q^{kn^2+n-in}a^n(1-x^iq^{2ni})(axq^{n+1})_\infty
(1/a)_n}{(q)_n(xq^n)_\infty}.\label{defH}
\end{eqnarray}
Andrews showed in \cite[p.\ 106-107]{ag} that
\begin{eqnarray}
H_{k,i}(a,x,q)-H_{k,i-1}(a,x,q)&=&x^{i-1}J_{k,k-i+1}(a,x,q)\label{H}\\
J_{k,i}(a,x,q)-J_{k,i-1}(a,x,q)&=&(xq)^{i-1}(J_{k,k-i+1}(a,xq,q)-aJ_{k,k-i+2}(a,xq,q)).
\end{eqnarray}

We plug $i=0$ in equation (\ref{defH}) and obtain $H_{k,0}=0$. 
We then plug $i=1$ in equation
\eqref{H} and obtain $H_{k,1}(a,x,q)=J_{k,k}(a,x,q)$. Then we plug $i=1$ 
in equation (\ref{defJ}) and
obtain
$J_{k,1}(a,x,q)=J_{k,k}(a,xq,q)$.
Finally we set $x=0$ in equations (\ref{defJ}) and (\ref{defH}) and get $J_{k,i}(a,0,q)=1$.
This implies for $1\le i\le k$, the following equations determine $J_{k,i}(a,x,q)$:
\begin{eqnarray*}
J_{k,i}(a,0,q)&=&1\\%  \ \ 1\le i\le k\\
J_{k,i}(a,x,q)-J_{k,i-1}(a,x,q)&=&(xq)^{i-1}(J_{k,k-i+1}(a,xq,q)-aJ_{k,k-i+2}(a,xq,q));\\
%&& \text{\hfill} 2\le i\le k\\
J_{k,1}(a,x,q)&=&J_{k,k}(a,xq,q).
\end{eqnarray*}

Let $\overline{\mathcal B}_{k,i}(a,x,q)=\sum_{n,j,p\ge 0} \overline{B}_{k,i}(n,j,p) a^j x^p q^n$ where $\overline{B}_{k,i}(n,j,p)$ is
the number of overpartitions counted by $\overline{B}_{k,i}(n,j)$ with $p$ parts.
We show that
\begin{lemma} For $1\le i\le k$,
$$
\overline{\mathcal B}_{k,i}(a,x,q)=J_{k,i}(-a,x,q).
$$
\end{lemma}

\begin{proof} %We give here a sketch. See \cite{lo} for a detailed proof.
The only overpartition with zero parts is the empty one.
Therefore $\overline{\mathcal B}_{k,i}(a,0,q)=1$. It is obvious that for $2\le i\le k$,
$\overline{\mathcal B}_{k,i}(a,x,q)-\overline{\mathcal B}_{k,i-1}(a,x,q)$ is the generating function for overpartitions such that $\forall \ell \text{, } f_{\ell} + f_{\ell+1}<k\overline{+1}$ and $f_1=i-1$.
Moreover $\overline{\mathcal B}_{k,i}(a,xq,q)$ is the generating function for overpartitions such that $\forall \ell \text{, } f_{\ell} + f_{\ell+1}<k\overline{+1}$, $f_2<i$ and $f_1=0$.
Therefore $(xq)^{i-1}\overline{\mathcal B}_{k,k-i+1}(a,xq,q)$ is the generating function for overpartitions such that $\forall \ell \text{, } f_{\ell} + f_{\ell+1}<k\overline{+1}$, $f_1=i-1$ and $1$ is not overlined and $a(xq)^{i-1}\overline{\mathcal B}_{k,k-i+2}(a,xq,q)$ is the generating function for overpartitions such that $\forall \ell \text{, } f_{\ell} + f_{\ell+1}<k\overline{+1}$, $f_1=i-1$ and $1$ is  overlined. We get, for $2\le i\le k$,
$\overline{\mathcal B}_{k,i}(a,x,q)-\overline{\mathcal B}_{k,i-1}(a,x,q)=(xq)^{i-1}(\overline{\mathcal B}_{k,k-i+1}(a,xq,q)+a\overline{\mathcal B}_{k,k-i+2}(a,xq,q)).$
Finally $\overline{\mathcal B}_{k,k}(a,xq,q)$ is the generating function for overpartitions such that $\forall \ell \text{, } f_{\ell} + f_{\ell+1}<k\overline{+1}$, $f_2<k$ and $f_1=0$
and therefore is equal to $\overline{\mathcal B}_{k,1}(a,x,q)$. \end{proof}

\begin{proof} We can now prove equation \eqref{mk}. The lemma implies that
\begin{equation*}
\overline{\mathcal B}_{k,i}(a,q)=J_{k,i}(-a,1,q).
\end{equation*}

Hence
\begin{eqnarray*}
\overline{\mathcal B}_{k,i}(a,q)&=&\frac{(-aq)_\infty}{(q)_\infty}\sum_{n=0}^{\infty} (-1)^n a^n \frac{q^{kn^2+n(k-i+1)}(-1/a)_n(1-q^{(2n+1)i})}{
(-aq)_{n+1}}\\
&&+aq\frac{(-aq)_\infty}{(q)_\infty}\sum_{n=0}^{\infty} (-1)^n a^n \frac{q^{kn^2+n(k-i+2)}(-1/a)_n(1-q^{(2n+1)(i-1)})}{
(-aq)_{n+1}}\\
&=&\frac{(-aq)_\infty}{(q)_\infty} \left(\sum_{n=0}^{\infty} (-1)^n a^n \frac{q^{kn^2+n(k+1)}(-1/a)_n (q^{-in}+aq^{1-(i-1)n})}{
(-aq)_{n+1}}\right.\\
&&-\left.\sum_{n=0}^{\infty} (-1)^n a^n \frac{q^{kn^2+n(k+1)}(-1/a)_n(q^{(n+1)i}+aq^{(n+1)(i-1)+1})}{
(-aq)_{n+1}}\right)\\
&=&\frac{(-aq)_\infty}{(q)_\infty}\left(\sum_{n=0}^{\infty} (-1)^n a^n \frac{q^{kn^2+n(k+1-i)}(-1/a)_n}{
(-aq)_n}\right.\\
&&-\left.\sum_{n=0}^{\infty} (-1)^n a^{n+1} \frac{q^{kn^2+n(k+i)+i}(-1/a)_{n+1}}{
(-aq)_{n+1}}\right)\\
&=&\frac{(-aq)_\infty}{(q)_\infty}\left(\sum_{n=0}^{\infty} (-1)^n a^n \frac{q^{kn^2+n(k+1-i)}(-1/a)_n}{
(-aq)_n}\right.\\
&&+\left.\sum_{n=-\infty}^{-1} (-1)^n a^{-n} \frac{q^{kn^2+n(k-i)}(-1/a)_{-n}}{
(-aq)_{-n}}\right)\\
&=&\frac{(-aq)_\infty}{(q)_\infty}\sum_{n=-\infty}^{\infty} (-1)^n a^n \frac{q^{kn^2+n(k+1-i)}(-1/a)_n}{
(-aq)_n}\\
&=&\overline{\mathcal E}_{k,i}(a,q)
\end{eqnarray*}
\end{proof}

\subsection{A combinatorial proof}\label{combproofmult}

This part is a generalization of \cite[Section 3]{bu1}. We only give a sketch of the proof
and details can be found in \cite{MAL}. In this section we represent
overpartitions by their multiplicity sequence $(f_0,f_1,f_2,\ldots )$. We include the multiplicity
$f_0$ to simplify the definitions, although it is always equal to 0.

We say that a sequence $(f_{m},\ldots, f_{m+\ell})$ is a {\em multuple} (tuple of multiplicity) if
\begin{itemize}
\item $f_{m+\ell}>0$,
\item $f_{m}$ is not overlined, and
\item $f_{m+p}$ is overlined for $1\le p\le \ell-1$.
\end{itemize}
The {\em length} of a multuple $(f_m,\ldots, f_{m+\ell})$ is
$\ell$ and its {\em weight} is $\sum_{i=0}^\ell (m+i)f_{m+i}$.
We divide a multiplicity sequence  of an overpartition
into multuples  going from the right to the left.
When we find a
positive multiplicity, we close a parenthesis to its right. We look for the next non-overlined
multiplicity to its left and open a parenthesis to the left of the multiplicity.
The {\em length} of a multiplicity sequence is the sum of the length of its multuples.
For partitions, the length is called
the number of pairs of multiplicities \cite{bu1}.
For example if the overpartition has for multiplicity sequence
$(0,\bar{2},0,2,\bar{1},1)$, then its multuples
are $((0,\bar{2}),0,(2,\bar{1},1))$. The first multuple has length 1 and the second 2. 
Therefore the length is $3$.\\

%We say that a
%multiplicity sequence $(f_m,\ldots ,f_{m+\ell})$
%is an {\em i-multuple} if and only if $m\ge 0$, $f_m>0$ and
%\begin{itemize}
%\item $f_m$ is not overlined and $\ell=1$ or
%\item $f_{m+p}$ is overlined for $0\le p\le \ell-1$ and $f_{m+\ell}=0$ or
%\item $f_{m+p}$ is overlined for $0\le p\le \ell-2$ and $f_{m+\ell-1}>0$ and is not overlined.
%\end{itemize}
%The {\em length} of an i-multuple $(f_m,\ldots, f_{m+\ell})$ is
%$\ell$ and its {\em weight} is $\sum_{i=0}^\ell (m+i)f_{m+i}$.
%We divide a multiplicity sequence  of an overpartition
%into i-multuples going from the left to the right.
%When we find a
%multiplicity $f_m>0$, we open a parenthesis to its left. If $f_m$ is
%not overlined then we close the parenthesis to the right of
%$f_{m+1}$. Otherwise, we look for the next non-overlined
%multiplicity, say $f_p$. If $f_p=0$ then we close the parenthesis to
%its right, otherwise we close the parenthesis to the right of
%$f_{p+1}$. Then we look for the next positive multiplicity, and so
%on. For example if the overpartition has for multiplicity sequence
%$(0,\bar{2},1,\bar{2},2,\bar{1},0,0,\bar{1},\bar{2},0)$, then the multuples
%are $(0,(\bar{2},1,\bar{2}),(2,\bar{1}),0,0,(\bar{1},\bar{2},0))$.

We define a map $F$ from  multuples of weight $n$ and length $\ell$
to  multiplicity sequences of weight $n-\ell$ and length $\ell$.
Given a multuple
$(f_m,\ldots, f_{m+\ell})$ then $F(f_m,\ldots, f_{m+\ell})$ is computed with the following algorithm :
\begin{itemize}
\item if $f_{m+\ell}=\overline{1}$, remove the overlining of $f_{m+\ell}$ and overline
$f_m$,
\item else if $\ell>1$, remove the overlining of $f_{m+\ell-1}$ and overline
$f_m$,
\item $f_{m+\ell} \leftarrow f_{m+\ell} -1$, and
\item $f_m \leftarrow f_m + 1$.
\end{itemize}
For example if $m=1$ and the multuple is $(1,\bar{1},\bar{3})$, we have $n=12$, $\ell=2$ and $F(1,\bar{1},\bar{3})=(\bar{2},1,\bar{2})$ whose weight is 10.\\

Let $\overline{B}_{k,i}(n,j,N)$ be the number of overpartitions
counted by $\overline{B}_{k,i}(n,j)$ of length $N$.
Let
$\overline{\mathcal
B}_{k,i}(N)=\sum_{n,j}\overline{B}_{k,i}(n,j,N)q^na^j$.
% Let $\overline{\beta}_{k,i-1}(N)$ (resp.\ $\overline{\delta}_{k,i}(N)$) be the generating functions for the %overpartitions
%obtained after applying $\beta^{-1}$
%(resp.  $\delta^{-1}$) to an overpartition counted by $\overline{\mathcal B}_{k,i}(N)$).
Now we divide the multiplicity sequence of an overpartition $\lambda\in \overline{\mathcal
B}_{k,i}(N)$ into multuples going from the
right to the left and we apply $F$ to each multuple. We call the result $F(\lambda)$.
If $F(\lambda)$ has a zero part, this part is discarded. Note that if $\lambda$ has weight $n$ and length $N$ then $F(\lambda)$ has weight $n-N$ and length $N$ or $N-1$.
For example,  $\lambda=(0,(0,\overline{1}),(1,\overline{1},\overline{3}))$ has weight  
24 and length 3 and 
$F(0,(0,\overline{1}),(1,\overline{1},\overline{3}))=(0,\overline{1},0,\overline{2},1,\overline{2})$
whose weight is 21 and length 3.

Let $\overline{\mathcal G}_{k,i}(N)$ be the generating function for the overpartitions
$\mu=F(\lambda)$ where $\lambda$ is an overpartition in $\overline{\mathcal B}_{k,{i+1}}(N)$
and has a multuple
$(f_0,\ldots,f_{\ell})$  with $\ell=1$ and $f_\ell\neq\bar{1}$.
\begin{proposition}
Let $\lambda$ be an overpartition in $\overline{\mathcal B}_{k,i}(N)$ with $N>0$. Then
\begin{itemize}
\item  $\lambda$ is an overpartition of $n$ and has a multuple
$(f_0,\ldots,f_{\ell})$ with $\ell>1$ or $f_\ell=\bar{1}$ if and only if
$F(\lambda)$ is an overpartition of $n-N$, has one less overlined part than $\lambda$ and is in  $\overline{\mathcal B}_{k,i}(N-1)$.
\item $\lambda$ is an overpartition of $n$  and has a multuple
$(f_0,\ldots,f_{\ell})$  with $\ell=1$ and $f_\ell\neq\bar{1}$ if and only if
$F(\lambda)$ is an overpartition of $n-N$, has the same number of overlined parts as $\lambda$ and
is in 
$\overline{\mathcal G}_{k,i-1}(N)$.
\item $\lambda$ is an overpartition of $n$  and has no multuple $(f_0,\ldots,f_{\ell})$ if and only if
$F(\lambda)$ is an overpartition of $n-N$, has the same number of overlined parts as $\lambda$ and
is in $\overline{\mathcal B}_{k,i+1}(N)$.
\end{itemize}

Let   $\lambda$ be an overpartition in $\overline{\mathcal
G}_{k,i}(N)$. Then
\begin{itemize}
\item $\lambda$ is an overpartition of $n$ and has a multuple
$(f_0,\ldots,f_{\ell})$ with $\ell>1$ or $f_\ell=\bar{1}$ if and only if
$F(\lambda)$ is an overpartition of $n-N$, has one less overlined part than $\lambda$ and is in  $\overline{\mathcal B}_{k,i}(N-1)$.
\item $\lambda$ is an overpartition of $n$  and has a multuple
$(f_0,\ldots,f_{\ell})$  with $\ell=1$ and $f_\ell\neq\bar{1}$  and
$F(\lambda)$ is an overpartition of $n-N$, has the same number of overlined parts as $\lambda$ and
is in $\overline{\mathcal G}_{k,i-1}(N)$.
\item $\lambda$ is an overpartition of $n$  and has no multuple $(f_0,\ldots,f_{\ell})$ if and only if
$F(\lambda)$ is an overpartition of $n-N+1$, has the same number of overlined parts as $\lambda$ and
is in $\overline{\mathcal B}_{k,i+1}(N-1)$.
\end{itemize}
\end{proposition}

\begin{proof} The full proof requires lots of details and  is given in \cite{MAL}.
We give  here the first ingredient of the proof, that is,
if $\lambda$ is an overpartition such that the condition
$\forall \ell$, $f_{\ell}+f_{\ell+1}<k\overline{+1}$ holds, then this condition still holds for $F(\lambda)$.
Indeed
the only successive multiplicities $f_m$ and $f_{m+1}$, where $f_m+f_{m+1}$ can increase
are such that $f_m$ or $f_{m+1}$ is the leftmost entry of a
multuple.  If  it is $f_m$ then after the operation
$f_m$ is now overlined and $f_{m}+f_{m+1}<k\overline{+1}$
still holds. If it is $f_{m+1}$ then
$f_{m}+f_{m+1}$ increases by one, only if $f_m=0$ before the operation.
We know that  $f_{m+1}<k-1$ before the operation (as $f_{m+1}+f_{m+2}<k$ and $f_{m+2}>0$)
then $f_{m}+f_{m+1}<k$ still holds after the operation.
\end{proof}

Note that $\overline{\mathcal B}_{k,i}(0)=1$ as the only overpartition of length 0 is the empty overpartition and that if an overpartition has a multuple $(f_0,\ldots ,f_\ell)$ then $i>1$. The previous remark and proposition imply that:
\begin{eqnarray*}
\overline{\mathcal B}_{k,1}(N) &=& q^N \overline{\mathcal B}_{k,2}(N); \\
\overline{\mathcal B}_{k,i}(N) &=& q^N (\overline{\mathcal B}_{k,i+1}(N) + \overline{\mathcal G}_{k,i-1}(N) +a\overline{\mathcal B}_{k,i}(N-1)), \text{ if $1<i\le k$};\\
\overline{\mathcal G}_{k,1}(N) &=& q^{N-1} \overline{\mathcal B}_{k,2}(N-1);\\
\overline{\mathcal G}_{k,i}(N) &=& q^{N-1} \overline{\mathcal B}_{k,i+1}(N-1) + q^N \overline{\mathcal G}_{k,i-1}(N) +aq^N \overline{\mathcal B}_{k,i}(N-1), \text{ if $1<i<k$};\\
\overline{\mathcal B}_{k,i}(0) &=& 1.\\
\end{eqnarray*}

Now we have all the ingredients to prove Proposition \ref{sec5}.
We remark that $\overline{\mathcal B}_{k,k+1}(N)=\overline{\mathcal B}_{k,k}(N) $ and we
set $\overline{\Gamma}_{k,i}(N)=\overline{\mathcal G}_{k,i}(N)+a\overline{\mathcal B}_{k,i+1}(N-1)$
if $i>0$ and $0$ otherwise. Then
\begin{eqnarray*}
\overline{\mathcal B}_{k,i}(N) &=& q^N (\overline{\mathcal B}_{k,i+1}(N) + \overline{\Gamma}_{k,i-1}(N)),  \text{ if $i<k$};\\
\overline{\mathcal B}_{k,k}(N) &=& \overline{\Gamma}_{k,k-1}(N)/(1-q^N);\\
\overline{\Gamma}_{k,i}(N) &=& (a+q^{N-1}) \overline{\mathcal B}_{k,i+1}(N-1) + q^N \overline{\Gamma}_{k,i-1}(N),  \text{ if $i>0$};\\
\overline{\mathcal B}_{k,i}(0) &=& 1;\\
\overline{\Gamma}_{k,0}(N) &=& 0.
\end{eqnarray*}
These are the recurrences proven in Proposition \ref{rec}.
Therefore
$
\overline{\mathcal B}_{k,i}(N)=\overline{\mathcal E}_{k,i}(N)
$
and Proposition \ref{sec5} is proved.

\section{Paths and successive Durfee squares}

We prove in this section that
\begin{equation}
\overline{\mathcal E}_{k,i}(a,q)=
\sum_{n_1 \ge \ldots \ge n_{k-1} \ge 0} \frac{q^{\binom{n_1+1}{2}+n_2^2+\cdots+n_{k-1}^2+n_i+\cdots+n_{k-1}}(-1/a)_{n_1}a^{n_1}}{(q)_{n_1-n_2}\cdots(q)_{n_{k-2}-n_{k-1}}(q)_{n_{k-1}}}.
\label{pathdur}
\end{equation}
This gives the case $\overline{E}_{k,i}(n,j)=\overline{D}_{k,i}(n,j)$ of Theorem \ref{main}
as the right hand side of Equation
\eqref{pathdur} is the generating function for the overpartitions counted by $\overline{D}_{k,i}(n,j)$
(see Equation \eqref{multseriesD}). We give an analytical and a combinatorial proof.

\subsection{An analytical proof}

We use
the Bailey lattice structure from \cite{Ag-An-Br1} to transform
 \eqref{pathdur} into  \eqref{gf}.  Recall that a
pair of sequences $(\alpha_n, \beta_n)$ form a Bailey pair with
respect to $a$ if for all $n \geq 0$ we have
$$
\beta_n = \sum_{r = 0}^n \frac{\alpha_r}{(q)_{n-r}(aq)_{n+r}}.
$$
We need the following lemma which is a special case of (3.8) in \cite{Ag-An-Br1}.:
\begin{lemma} \label{lattice}
If $(\alpha_n,\beta_n)$ is a Bailey pair with respect to $q$, then
for all $0 \leq i \leq k$ we have
\begin{eqnarray}
\frac{1}{(q,-aq)_{\infty}} &\times& \sum_{n_1
\geq \cdots \geq n_k \geq 0} \frac{q^{\binom{n_1+1}{2} + n_2^2 + \cdots n_k^2 +
n_{i+1} + \cdots
n_k}(-1/a)_{n_1}a^{n_1}}{(q)_{n_1-n_2}\cdots(q)_{n_{k-1} -
n_k}} \beta_{n_k} \nonumber \\
&=& \frac{\alpha_0}{(q)_{\infty}^2} +
\frac{1}{(q)_{\infty}^2}\sum_{n \geq 1}
\frac{(-1/a)_na^nq^{(n^2-n)(i-1/2)+in}(1-q)}{(-aq)_n} \nonumber \\
&\times& \left(\frac{q^{(n^2+n)(k-i)}}{(1-q^{2n+1})}\alpha_n -
\frac{q^{((n-1)^2 + (n-1))(k-i)+ 2n-1}}{(1-q^{2n-1})} \alpha_{n-1}
\right)\label{lattic}
\end{eqnarray}
\end{lemma}

\begin{proof}
We set $a = q$, $\rho_1 = -1/a$,
and then let $n$ as well as all remaining $\rho_i$ and $\sigma_i$
tend to $\infty$ in (3.8) of \cite{Ag-An-Br1} to obtain \eqref{lattic}.
\end{proof}

\noindent \emph{Proof of \eqref{pathdur}.}  We use the Bailey pair with respect to $q$
\cite[p.468, (B3)]{Sl1},
$$
\beta_n = \frac{1}{(q)_{\infty}} \hskip.5in \text{and} \hskip.5in
\alpha_n = \frac{(-1)^nq^{n(3n+1)/2}(1-q^{2n+1})}{(1-q)}.
$$
Substituting into Lemma \ref{lattice} and simplifying, we obtain
\begin{align*}
\sum_{n_1 \ge \cdots \ge n_{k} \ge 0} & \frac{q^{\binom{n_1+1}{2} + n_2^2 +
\cdots + n_{k}^2 + n_{i+1} + \cdots +
n_{k}}(-1/a)_{n_{1}}a^{n_1}}{(q)_{n_1-n_2}\cdots(q)_{n_{k-1}
- n_{k}} (q)_{n_k}} \\ =&
\frac{(-aq)_{\infty}}{(q)_{\infty}} \left(1 + \sum_{n =
1} ^\infty \frac{q^{(k+1)n^2+(k-i+1)n +
}(-a)^n(-1/a)_n}{(-aq)_n} \right. \\  +& \left.
\sum_{n = 1} ^\infty \frac{q^{(k+1)n^2 -(k-i)n}(-a)^n(-1/a)_n}{(-aq)_n} \right)\\
=&
\frac{(-aq)_{\infty}}{(q)_{\infty}} \left(1 + \sum_{n=1}^\infty  \frac{q^{(k+1)n^2+(k-i+1)n
}(-a)^n(-1/a)_n}{(-aq)_n} \right. \\  +& \left.
\sum_{n=-\infty}^{-1} \frac{q^{(k+1)n^2 +(k-i+1)n}(-a)^n(-1/a)_n}{(-aq)_n} \right).
\end{align*}
Replacing $k$ by $k-1$ and $i$ by $i-1$ gives \eqref{gf}.
\qed

\subsection{A combinatorial proof}

We will use and generalize the notion of  relative height of a peak. This notion was defined by Bressoud in \cite{b}
for paths with no South steps
and a simpler version of the definition was given in \cite{bep}.
We adapt this definition for the paths with South steps.
\begin{definition}
The relative height of a peak $(x,y)$ is the largest integer $h$ for which we can find two vertices on the path, $(x',y-h)$ and $(x'',y-h)$, such that $x'<x\le x''$ and such that between these two vertices there are no peaks of height $> y$ and every peak of height $y$ has abscissa $\ge x$.
\end{definition}
The original definition was the same except $x\le x''$ was $x<x''$. Indeed when there are no South steps
the case $x=x''$ is impossible.

%As an example, the relative heights in the path from Figure \ref{exrhp} are, from left to right: 1, 4, 1, 2, 4.
%\begin{figure}[ht]
%\centering
%\psset{unit=0.4}
%\begin{pspicture}(0,0)(26,4)
%\psaxes[labels=none,ticksize=1pt](26,4)
%\psline(0,2)(1,3)(4,0)(8,4)(9,3)(10,4)(14,0)(17,3)(19,1)(22,4)(26,0)
%\psdots[dotstyle=x](2,2)
%\psdots[dotstyle=x](3,1)
%\psdots[dotstyle=x](5,1)
%\psdots[dotstyle=x](6,2)
%\psdots[dotstyle=x](7,3)
%\psdots[dotstyle=x](11,3)
%\psdots[dotstyle=x](12,2)
%\psdots[dotstyle=x](13,1)
%\psdots[dotstyle=x](15,1)
%\psdots[dotstyle=x](16,2)
%\psdots[dotstyle=x](18,2)
%\psdots[dotstyle=x](20,2)
%\psdots[dotstyle=x](21,3)
%\psdots[dotstyle=x](23,3)
%\psdots[dotstyle=x](24,2)
%\psdots[dotstyle=x](25,1)
%\rput(1,3.5){\small 1}
%\rput(8,4.5){\small 4}
%\rput(10,4.5){\small 1}
%\rput(17,3.5){\small 2}
%\rput(22,4.5){\small 4}
%\end{pspicture}
%\caption{Example of a path and its relative heights.}
%\label{exrhp}
%\end{figure}

%For the paths corresponding to overpartitions, that is to say the paths counted by $\overline{E}_{k,i}(n,j)$, we have to modify the definition of the relative height a little bit to take into account the NES peaks, for which we can have $x''=x$ (this is the case for the first peak of the path in Figure \ref{exrho}, for example).

%\begin{definition}
%The relative height of a peak $(x,y)$ is the largest integer $h$ for which we can find two vertices on the path, $(x',y-h)$ and $(x'',y-h)$, such that $x' < x \le x''$ and such that between these two vertices there are no peaks of height $> y$ and every peak of height $y$ has abscissa $\ge x$.
%\end{definition}

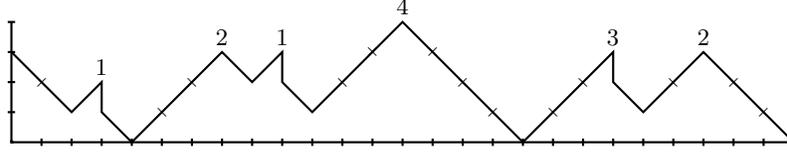
\begin{figure}[ht]
\centering
\psset{unit=0.4}
\begin{pspicture}(0,0)(26,4)
\psaxes[labels=none,ticksize=1pt](26,4)
\psline(0,3)(2,1)(3,2)(3,1)(4,0)(7,3)(8,2)(9,3)(9,2)(10,1)(13,4)(17,0)(20,3)(20,2)(21,1)(23,3)(26,0)
\psdots[dotstyle=x](1,2)
\psdots[dotstyle=x](5,1)
\psdots[dotstyle=x](6,2)
\psdots[dotstyle=x](11,2)
\psdots[dotstyle=x](12,3)
\psdots[dotstyle=x](14,3)
\psdots[dotstyle=x](15,2)
\psdots[dotstyle=x](16,1)
\psdots[dotstyle=x](18,1)
\psdots[dotstyle=x](19,2)
\psdots[dotstyle=x](22,2)
\psdots[dotstyle=x](24,2)
\psdots[dotstyle=x](25,1)
\rput(3,2.5){\small 1}
\rput(7,3.5){\small 2}
\rput(9,3.5){\small 1}
\rput(13,4.5){\small 4}
\rput(20,3.5){\small 3}
\rput(23,3.5){\small 2}
\end{pspicture}
\caption{An example of a path and its relative heights.}
\label{exrho}
\end{figure}

For example, the relative heights in the path from Figure 
\ref{exrho} are, from left to right: 1, 2, 1, 4, 3, 2.

We will prove here that
\begin{proposition}\label{multseries}
The coefficient of $q^na^\ell$ in
\begin{equation}\label{ims}
\frac{q^{\binom{n_1+1}{2}+n_2^2+\cdots+n_{k-1}^2+n_i+\cdots+n_{k-1}}(-1/a)_{n_1}a^{n_1}}{(q)_{n_1-n_2}\cdots(q)_{n_{k-2}-n_{k-1}}(q)_{n_{k-1}}}
\end{equation}
is the number of paths with major index $n$ and $\ell$ South steps, starting at height $k-i$, whose height is less than 
$k$ and having $n_j$ peaks of relative height $\ge j$ for $1 \le j \le k-1$. Therefore,
$$
\sum_{n_1 \ge \cdots \ge n_{k-1} \ge 0} \frac{q^{\binom{n_1+1}{2}+n_2^2+\cdots+n_{k-1}^2+n_i+\cdots+n_{k-1}}(-1/a)_{n_1}a^{n_1}}{(q)_{n_1-n_2}\cdots(q)_{n_{k-2}-n_{k-1}}(q)_{n_{k-1}}} = \overline{\mathcal E}_{k,i}(a,q).
$$
\end{proposition}

To prove this proposition we will use a result of Bressoud \cite{b}
\begin{lemma}\label{lembressoud}
The coefficient of $q^n$ in
$$
\frac{q^{n_1^2+n_2^2+\cdots+n_{k-1}^2+n_i+\cdots+n_{k-1}}}{(q)_{n_1-n_2}\cdots(q)_{n_{k-2}-n_{k-1}}(q)_{n_{k-1}}}
$$
is the number of paths with major index $n$, no South steps, starting at height $k-i$, 
whose height is less than $k$ and having $n_j$ peaks of relative height $\ge j$ for $1 \le j \le k-1$.
\end{lemma}

\begin{figure}[ht]
        \psset{unit=0.4}
        \begin{pspicture}(0,0)(13,3)
        \psaxes[labels=none,ticksize=1pt](13,3)
        \psline(0,3)(2,1)(3,2)(5,0)(8,3)(9,2)(10,3)(13,0)
        \psdots[dotstyle=x](1,2)
        \psdots[dotstyle=x](4,1)
        \psdots[dotstyle=x](6,1)
        \psdots[dotstyle=x](7,2)
        \psdots[dotstyle=x](11,2)
        \psdots[dotstyle=x](12,1)
        \rput(3,2.5){\small 1}
        \rput(8,3.5){\small 3}
        \rput(10,3.5){\small 1}
        \end{pspicture}
        \caption{Another example of a path and its relative heights.}
        \label{expathBressoud}
\end{figure}
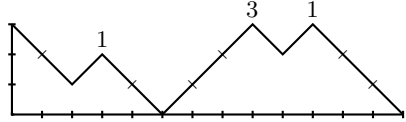

An example of such a path, taken from \cite{b}, is shown on Figure \ref{expathBressoud}. For that path, we have $k=4$, $i=1$, $n_1=3$, $n_2=1$ and $n_3=1$.

We can now move on to the proof of Proposition \ref{multseries}.
\begin{proof}
We generalize the argument of Bressoud in \cite{b}.
Consider a path with no South steps that starts at height $k-i$, whose height is less than $k-1$ and that has $n_j$ peaks of relative height $\ge j-1$ for $2 \le j \le k-1$.
By Lemma \ref{lembressoud}, such paths are counted by
\begin{equation}\label{bih}
\mathcal{P}_{k,i}(q)=\frac{q^{n_2^2+\cdots+n_{k-1}^2+n_i+\cdots+n_{k-1}}}{(q)_{n_2-n_3}\cdots(q)_{n_{k-2}-n_{k-1}}(q)_{n_{k-1}}}
\end{equation}
where $2 \le i \le k$.
Note that we have replaced $n_j$ by $n_{j+1}$ so that the terms in $n_1$ should be introduced by the algorithm described below, which will then give us the generating function in the same form as in \eqref{ims}.

% les entrees et la sortie de l'algo
For any given $k \ge 2$ and $i$ such that $1 \le i \le k$, we describe an algorithm which generates a path counted by \eqref{ims} from a path counted by $\mathcal{P}_{k,i}(q)$ if $i \ge 2$ or $\mathcal{P}_{k,2}(q)$ if $i=1$, a partition $\lambda$ into distinct parts which lie in $[0,n_1-1]$ and a partition $b$ into $n_1-n_2$ nonnegative parts.
%We describe an algorithm which generates a path counted by \eqref{ims} from a path counted by \eqref{bih}, a partition $\lambda$ into distinct parts which lie in $[0,n_1-1]$ and a partition $b$ into $n_1-n_2$ nonnegative parts.

% ce qu'il faut prouver
We will need to prove that this contruction is uniquely reversible, that the algorithm generates all of our paths, that the distribution of relative heights is not modified (except at the first step where all the peaks will be raised by one) and that the algorithm affects the generating function in the appropriate way.

% l'algo
We first perform a ``volcanic uplift'' by inserting  a NES peak at each peak (see Figure \ref{volcanic}).
This increases all the relative height by one.

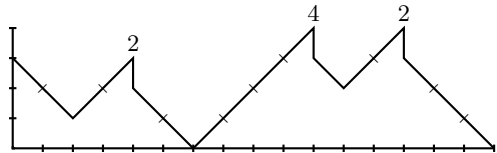
\begin{figure}[ht]
        \centering
        \psset{unit=0.4}
        \begin{pspicture}(0,0)(16,4)
        \psaxes[labels=none,ticksize=1pt](16,4)
        \psline(0,3)(2,1)(4,3)(4,2)(6,0)(10,4)(10,3)(11,2)(13,4)(13,3)(16,0)
        \psdots[dotstyle=x](1,2)
        \psdots[dotstyle=x](3,2)
        \psdots[dotstyle=x](5,1)
        \psdots[dotstyle=x](7,1)
        \psdots[dotstyle=x](8,2)
        \psdots[dotstyle=x](9,3)
        \psdots[dotstyle=x](12,3)
        \psdots[dotstyle=x](14,2)
        \psdots[dotstyle=x](15,1)
        \rput(4,3.5){\small 2}
        \rput(10,4.5){\small 4}
        \rput(13,4.5){\small 2} 
        \end{pspicture}
        \caption{Effect of the volcanic uplift on the path from Figure \ref{expathBressoud}.}
        \label{volcanic}
\end{figure}

We then insert $n_1-n_2$ NES peaks at the beginning of the path (see Figure \ref{insertpeaks}). Note that all these peaks have relative height one and that they are the only peaks of relative height one since the volcanic uplift has increased all the relative heights by one.

\begin{figure}[ht]
        \psset{unit=0.4}
        \begin{pspicture}(0,0)(20,4)
        \psaxes[labels=none,ticksize=1pt](20,5)
        \psline(0,3)(1,4)(1,3)(2,4)(2,3)(3,4)(3,3)(4,4)(4,3)(6,1)(8,3)(8,2)(10,0)(14,4)(14,3)(15,2)(17,4)(17,3)(20,0)
        \psdots[dotstyle=x](5,2)
        \psdots[dotstyle=x](7,2)
        \psdots[dotstyle=x](9,1)
        \psdots[dotstyle=x](11,1)
        \psdots[dotstyle=x](12,2)
        \psdots[dotstyle=x](13,3)
        \psdots[dotstyle=x](16,3)
        \psdots[dotstyle=x](18,2)
        \psdots[dotstyle=x](19,1)
        \rput(1,4.5){\small 1}
        \rput(2,4.5){\small 1}
        \rput(3,4.5){\small 1}
        \rput(4,4.5){\small 1}
        \rput(8,3.5){\small 2}
        \rput(14,4.5){\small 4}
        \rput(17,4.5){\small 2}         
        \end{pspicture}
        \caption{After adding the $n_1-n_2=4$ NES peaks of relative height one to the path from Figure \ref{volcanic}.}
        \label{insertpeaks}
\end{figure}
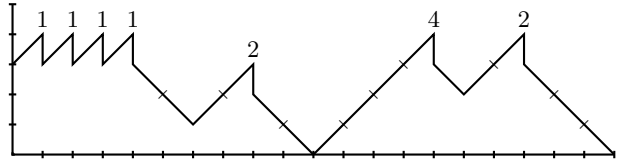

If $i=1$, we introduce an extra SE step at the beginning of the path, from $(0,k-1)$ to $(1,k-2)$.

Now if $\lambda$ contains a part $j-1$ ($1 \leq j \leq n_1$), we transform the $j$th NES peak from the right into a NESE peak (see Figure \ref{transfpeaks}).

\begin{figure}[ht]
        \psset{unit=0.4}
        \begin{pspicture}(0,0)(24,4)
        \psaxes[labels=none,ticksize=1pt](24,4)
        \psline(0,3)(1,4)(1,3)(2,4)(3,3)(4,4)(5,3)(6,4)(9,1)(11,3)(11,2)(13,0)(17,4)(19,2)(21,4)(21,3)(24,0)
        \psdots[dotstyle=x](7,3)
        \psdots[dotstyle=x](8,2)
        \psdots[dotstyle=x](10,2)
        \psdots[dotstyle=x](12,1)
        \psdots[dotstyle=x](14,1)
        \psdots[dotstyle=x](15,2)
        \psdots[dotstyle=x](16,3)
        \psdots[dotstyle=x](18,3)
        \psdots[dotstyle=x](20,3)
    \psdots[dotstyle=x](22,2)
    \psdots[dotstyle=x](23,1)
        \rput(1,4.5){\small 1}
        \rput(2,4.5){\small 1}
        \rput(4,4.5){\small 1}
        \rput(6,4.5){\small 1}
        \rput(11,3.5){\small 2}
        \rput(17,4.5){\small 4}
        \rput(21,4.5){\small 2}         
         \end{pspicture}
        \caption{Effect of transforming some NES peaks into NESE peaks in the path from Figure \ref{insertpeaks}. In this example, $\lambda=(5,4,3,1)$.}
        \label{transfpeaks}
\end{figure}
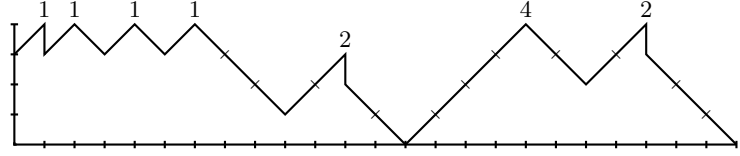

Finally, for $1 \leq j \leq n_1-n_2$, we move the $j$th peak of relative height one from the right $b_j$ times according to the rules illustrated in Figure \ref{rules}.

\begin{figure}[ht]
        \psset{unit=0.4}
        \begin{pspicture}(0,-13)(10,1)

        \psline(0,0)(1,1)(2,0)(3,0)
        \psline(7,0)(8,0)(9,1)(10,0)
        \psline[arrows=->](4,0.5)(6,0.5)
        \rput[b](5,1){1}

        \rput(-1,-3){
        \psline(0,0)(1,1)(2,0)(4,2)
        \psdots[dotstyle=x](3,1)
        \psline(8,0)(10,2)(11,1)(12,2)
        \psdots[dotstyle=x](9,1)
        \psline[arrows=->](5,1)(7,1)
        \rput[b](6,1.5){2}}

        \rput(-1,-5.5){
        \psline(0,2)(1,1)(2,2)(4,0)
        \psdots[dotstyle=x](3,1)
        \psline(8,2)(10,0)(11,1)(12,0)
        \psdots[dotstyle=x](9,1)
        \psline[arrows=->](5,1)(7,1)
        \rput[b](6,1.5){3}}

        \rput(1,-7.5){
        \psline(0,0)(1,1)(1,0)(2,0)
        \psline(6,0)(7,0)(8,1)(8,0)
        \psline[arrows=->](3,0.5)(5,0.5)
        \rput[b](4,1){4}}

        \rput(0,-10.5){
        \psline(0,0)(1,1)(1,0)(3,2)
        \psdots[dotstyle=x](2,1)
        \psline(7,0)(9,2)(9,1)(10,2)
        \psdots[dotstyle=x](8,1)
        \psline[arrows=->](4,1)(6,1)
        \rput[b](5,1.5){5}}

        \rput(0,-13){
        \psline(0,2)(1,1)(2,2)(2,1)(3,0)
        \psline(7,2)(9,0)(10,1)(10,0)
        \psdots[dotstyle=x](8,1)
        \psline[arrows=->](4,1)(6,1)
        \rput[b](5,1.5){6}}

        \end{pspicture}
\caption{The rules for moving peaks of relative height one.}
\label{rules}
\end{figure}
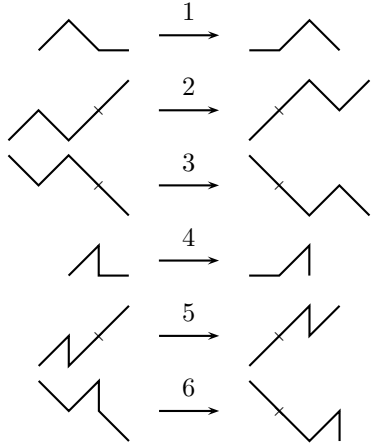

When we move a peak, it can meet the next peak to the right. We say that a peak $(x,y)$ meets a peak $(x',y')$ if
$$
x' - x =
\begin{cases}
2 \text{ if $(x,y)$ is a NESE peak}\\
1 \text{ if $(x,y)$ is a NES peak}
\end{cases}.
$$
If this happens, we abandon the peak we have been moving and move the next peak to the right.
If we come up against a sequence of adjacent peaks, we move the rightmost peak in the sequence (see Figure \ref{collision}).

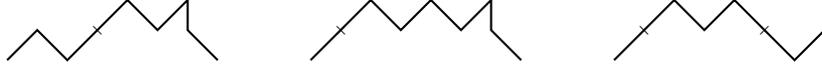
\begin{figure}
        \centering
        \psset{unit=0.4}
        \begin{pspicture}(0,0)(7,2)
        \psline(0,0)(1,1)(2,0)(4,2)(5,1)(6,2)(6,1)(7,0)
        \psdots[dotstyle=x](3,1)
        \end{pspicture}
        \hspace{1cm}
        \begin{pspicture}(0,0)(7,2)
        \psline(0,0)(2,2)(3,1)(4,2)(5,1)(6,2)(6,1)(7,0)
        \psdots[dotstyle=x](1,1)
        \end{pspicture}
        \hspace{1cm}
        \begin{pspicture}(0,0)(7,2)
        \psline(0,0)(2,2)(3,1)(4,2)(6,0)(7,1)(7,0)
        \psdots[dotstyle=x](1,1)
        \psdots[dotstyle=x](5,1)
        \end{pspicture}
        \caption{We want to move the leftmost peak to the right twice, but after the first move, we come up against a sequence of adjacent peaks. We then move the rightmost peak in this sequence.}
        \label{collision}
\end{figure}
To conclude the proof, we must show that the distribution of relative heights is not modified by the operations of Figure \ref{rules} (Proposition \ref{distrh}), that the construction procedure is uniquely reversible (Proposition \ref{inv}), that we generate all of our paths (Proposition \ref{genallpaths}), and that the algorithm affects the generating function in the appropriate way (Proposition \ref{gfaffect}).
\end{proof}

\begin{proposition}\label{distrh}
The operations of Figure \ref{rules} preserve the number of peaks of relative height $\geq j$ for all $j$.
\end{proposition}
\begin{proof}
Let us show it for each operation. We call $p$ the peak which is moved. Remember that before the move, the relative height of $p$ is 1.

For operations 1 and 4, the relative height of $p$ clearly remains 1 after the move.
The other peaks are not affected and their relative heights are therefore not modified.

For operations 2 and 5, it can be easily shown that the relative height of $p$ remains 1 (see Figure \ref{op25} for an example) unless it meets a peak, in which case the two peaks will swap their relative heights (see Figure \ref{invrh} for an example). In both cases, the relative heights of the other peaks are not modified.

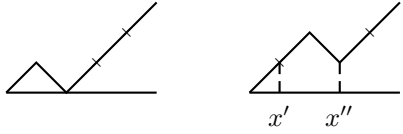
\begin{figure}[ht]
\centering
\psset{unit=0.4}
\begin{pspicture}(0,-1)(5,3)
\psline(0,0)(5,0)
\psline(0,0)(1,1)(2,0)(5,3)
\psdots[dotstyle=x](3,1)
\psdots[dotstyle=x](4,2)
\end{pspicture}
\hspace{1cm}
\begin{pspicture}(0,-1)(5,3)
\psline(0,0)(5,0)
\psline(0,0)(2,2)(3,1)(5,3)
\psdots[dotstyle=x](1,1)
\psdots[dotstyle=x](4,2)
\psline[linestyle=dashed](1,1)(1,0)
\psline[linestyle=dashed](3,1)(3,0)
\uput[d](1,0){$x'$}
\uput[d](3,0){$x''$}
\end{pspicture}
\caption{Case where the relative height of $p$ is not modified.}
\label{op25}
\end{figure}

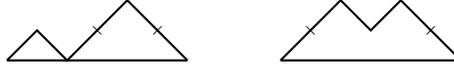
\begin{figure}[ht]
\centering
\psset{unit=0.4}
\begin{pspicture}(0,0)(6,2)
\psline(0,0)(6,0)
\psline(0,0)(1,1)(2,0)(4,2)(6,0)
\psdots[dotstyle=x](3,1)
\psdots[dotstyle=x](5,1)
\end{pspicture}
\hspace{1cm}
\begin{pspicture}(0,0)(6,2)
\psline(0,0)(6,0)
\psline(0,0)(2,2)(3,1)(4,2)(6,0)
\psdots[dotstyle=x](1,1)
\psdots[dotstyle=x](5,1)
\end{pspicture}
\caption{When $p$ meets a peak, their relative heights are inverted.
This does not modify the number of peaks with a given relative height.}
\label{invrh}
\end{figure}

For operations 3 and 6, the relative height of $p$ will clearly remain 1 since the peak located to its left has a larger (or equal) height (see Figure \ref{op36}).

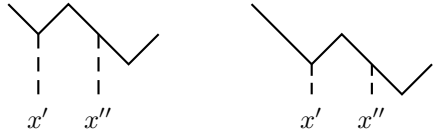
\begin{figure}[ht]
\psset{unit=0.4}
\begin{pspicture}(0,0)(5,-4)
\psline(0,0)(1,-1)(2,0)(4,-2)(5,-1)
\psline[linestyle=dashed](1,-1)(1,-3)
\psline[linestyle=dashed](3,-1)(3,-3)
\uput[d](1,-3){$x'$}
\uput[d](3,-3){$x''$}
\end{pspicture}
\hspace{1cm}
\begin{pspicture}(0,0)(6,-4)
\psline(0,0)(2,-2)(3,-1)(5,-3)(6,-2)
\psline[linestyle=dashed](2,-2)(2,-3)
\psline[linestyle=dashed](4,-2)(4,-3)
\uput[d](2,-3){$x'$}
\uput[d](4,-3){$x''$}
\end{pspicture}
\caption{The relative height of $p$ remains one when we apply the operation 3 or 6.}
\label{op36}
\end{figure}
\end{proof}

\begin{proposition}\label{inv}
The construction used in the proof of Proposition \ref{multseries} is uniquely reversible.
\end{proposition}
\begin{proof}
We describe the algorithm to undo the algorithm used to prove Proposition \ref{multseries}.
We start with a path counted
by \eqref{ims}.

We first have to move the peaks of relative height 1 to the left. We
begin with the leftmost of these peaks and we move it to the left
end of the path. If, during this move, our peak becomes adjacent to
another peak, we abandon the peak we have been moving and we move
the peak to its left (we do the same if we come up against a
sequence of contiguous peaks). The number of moves we had to perform
to bring the peak to the beginning of the path gives us $b_1$. We
proceed similarly for the other peaks of relative height 1, which
gives us a partition $b$ into $n_1-n_2$ parts $\geq 0$.
Note that we can do this because of Proposition \ref{distrh}: since the number of peaks of relative height 1 is preserved when we move the peaks, the peaks of relative height 1 in the path counted by \eqref{ims} correspond to these which were added by the volcanic uplift and the subsequent insertion of peaks.

At this step, among the $n_1$ peaks, some are NESE peaks. We transform these
NESE peaks into NES peaks: if the $j$th peak from the right was
transformed, it gives a part $j-1$ in $\lambda$.

If $i=1$, we remove the first South-East step of the path. Finally, we remove
the $n_1-n_2$ NES peaks at the beginning of the path and we decrease by 1 the height of the remaining peaks by removing the NES
peaks. The resulting path is counted by $\mathcal{P}_{k,2}(q)$ if $i=1$ and by $\mathcal{P}_{k,i}(q)$ otherwise.
\end{proof}

\begin{proposition}\label{genallpaths}
Any path in $\overline{E}_{k,i}(n,j)$ can be generated by our algorithm.
\end{proposition}
\begin{proof}
This is easy to see using the reverse algorithm described in the proof of Proposition \ref{inv}. That algorithm can be applied to any path $P$ counted by \eqref{ims} and gives a Bressoud path $P'$ (i.e.\ a path counted by \eqref{bih}). If we apply the direct algorithm to $P'$, we will obtain our initial path $P$ counted by \eqref{ims}. Thus, for any path $P$ in $\overline{E}_{k,i}(n,j)$, there exists a Bressoud path $P'$ which gives that path.
\end{proof}

\begin{proposition}\label{gfaffect}
Our algorithm affects the generating function in the appropriate way.
\end{proposition}
\begin{proof}
The volcanic uplift increases the major index of the path by
$$
1 + 2 + \cdots + n_2 = \binom{n_2+1}{2}
$$
and the relative height of each peak by one. Moreover, the $n_2$ NES peaks introduce a factor $a^{n_2}$.

The new peaks introduced after the uplift have total major index $\binom{n_1-n_2+1}{2}$ and they increase the abscissa of each of the old peaks by $n_1-n_2$. Since they are NES peaks, they also give a factor $a^{n_1-n_2}$. Altogether, the two operations introduce a factor
$$
q^{\binom{n_2+1}{2}} a^{n_2} \times q^{\binom{n_1-n_2+1}{2} + n_2(n_1-n_2)} a^{n_1-n_2} = q^{\binom{n_1+1}{2}} a^{n_1}.
$$

If $i=1$, we add an extra SE step at the beginning of the path, which introduces a factor $q^{n_1}$.

Transforming the $j$th peak from the right into a NESE peak increases the major index of the path by $j-1$ because the $j-1$ rightmost peaks are shifted by 1 to the right. We do that if there is a part $j-1$ in $\lambda$; altogether, the major index of the path increases precisely by the size of $\lambda$, which is a partition into distinct parts in $[0,n_1-1]$. Since we transform a NES peak into a NESE peak for each part of $\lambda$, this step introduces a factor $(-1/a)_{n_1}$.

Finally, when we move the $j$th peak of relative height one from the right $b_j$ times, we increase the abscissa of the path by $b_j$. Altogether, the major index  of the path increases precisely by the size of $b=(b_1,b_2,\ldots,b_{n_1-n_2})$, which is a partition into $n_1-n_2$ nonnegative parts. Such partitions are counted by $\frac{1}{(q)_{n_1-n_2}}$.
\end{proof}

The multiple series in Proposition \ref{multseries} can be
re-expressed as \eqref{multseriesD}, which is the generating
function for overpartitions with $i-1$ successive Durfee squares
followed by $k-i$ successive Durfee rectangles, the first one
being a generalized Durfee square/rectangle. %For $i=k$, this is
%Theorem \ref{start}. {\bf FIX THIS}.

\section{New partition theorems}

We first prove Corollary \ref{aftergf} and then extract its combinatorial information.

\subsection{Proof of Corollary \ref{aftergf}}

We recall here the Jacobi Triple Product identity (JTP) given in equation \eqref{JTP}:
\begin{equation*}
(-1/z,-zq,q;q)_\infty=\sum_{n=-\infty}^{\infty}z^nq^{\binom{n+1}{2}}.
\end{equation*}
and the result of Theorem \ref{thegf}
$$
\overline{\mathcal E}_{k,i}(a,q)=\frac{(-aq)_\infty}{(q)_\infty}\sum_{n=-\infty}^{\infty}(-1)^n a^n
q^{2k\binom{n+1}{2}-n(i-1)}\frac{(-1/a)_n}{(-aq)_n}.
$$

We first prove \eqref{andrews}.  Using Theorem \ref{thegf}, we get
$$
\overline{\mathcal E}_{k,i}(0,q)=\frac{1}{(q)_\infty}\sum_{n=-\infty}^{\infty}
(-1)^n q^{(2k+1)\binom{n+1}{2}-ni}.
$$
We substitute $q\rightarrow q^{2k+1}$, $z\rightarrow -q^{-i}$ in \eqref{JTP} and
get
$$
\overline{\mathcal E}_{k,i}(0,q)=\frac{(q^i,q^{2k+1-i},q^{2k+1};q^{2k+1})_\infty}{(q)_\infty}.
$$

Second we prove  \eqref{gordon}.  Using Theorem \ref{thegf}, we get
$$
\overline{\mathcal E}_{k,i}(1/q,q^2)=\frac{(-q;q^2)_\infty}{(q^2;q^2)_\infty}
\sum_{n=-\infty}^{\infty}(-1)^n q^{4k\binom{n+1}{2}-n(2i-1)}.
$$
We substitute $q\rightarrow q^{4k}$, $z\rightarrow -q^{-2i+1}$ in \eqref{JTP} and
get
$$
\overline{\mathcal E}_{k,i}(1/q,q^2)=\frac{(q^2;q^4)_\infty(q^{2i-1},q^{4k+1-2i},q^{4k};q^{4k})_\infty}
{(q)_\infty}.
$$

Third we prove \eqref{gord}.
Using Theorem \ref{thegf}, we get
\begin{eqnarray*}
\overline{\mathcal E}_{k,i}(1,q)&=&\frac{(-q)_\infty}{(q)_\infty}\sum_{n=-\infty}^{\infty}
(-1)^n q^{2k\binom{n+1}{2}-n(i-1)}\frac{2}{1+q^n}\\
\end{eqnarray*}
Writing half of this series twice, once with $n$ and once with $-n$, we have 
\begin{eqnarray*}
\overline{\mathcal E}_{k,i}(1,q)&=&\frac{(-q)_\infty}{(q)_\infty}
%(1+\sum_{n=1}^\infty
%(-1)^n q^{2k{n+1\choose 2}-ni)}\frac{2q^n}{1+q^n}
%+(-1)^n q^{2k{n+1\choose 2}+n(i-2k)}\frac{2}{1+q^{n}}\\
{\small \left(\sum_{n=-\infty}^{\infty}(-1)^n q^{2k\binom{n+1}{2}-ni}\frac{q^n}{1+q^n}
+\sum_{n=-\infty}^{\infty}(-1)^n q^{2k\binom{n+1}{2}+n(i-2k)}\frac{1}{1+q^{n}}\right)}
\end{eqnarray*}
Therefore
$$
\overline{\mathcal E}_{k,k}(1,q)=\frac{(-q)_\infty}{(q)_\infty}\sum_{n=-\infty}^{\infty}(-1)^n q^{2k\binom{n+1}{2}-nk};
$$
and for $i<k$
\begin{eqnarray*}
\overline{\mathcal E}_{k,i}(1,q)&=&
\frac{(-q)_\infty}{(q)_\infty}
{\small 
\left(\sum_{n=-\infty}^{\infty}(-1)^n q^{2k\binom{n+1}{2}-ni}
-\sum_{n=-\infty}^{\infty}(-1)^n q^{2k\binom{n+1}{2}-n(i+1)}\frac{q^n}{1+q^n}\right.}\\
&&{\small \left.
+\sum_{n=-\infty}^{\infty}(-1)^n q^{2k\binom{n+1}{2}+n(i-2k)}
-\sum_{n=-\infty}^{\infty}(-1)^n q^{2k\binom{n+1}{2}+n(i+1-2k)}\frac{1}{1+q^{n}}
\right)}\\
&=&
\frac{(-q)_\infty}{(q)_\infty}
\left(\sum_{n=-\infty}^{\infty}(-1)^n q^{2k\binom{n+1}{2}}(q^{-ni}+q^{n(i-2k)})
\right)-\overline{\mathcal E}_{k,i+1}(1,q)\\
&=&\frac{(-q)_\infty}{(q)_\infty}\sum_{j=0}^{2(k-i)}(-1)^j\sum_{n=-\infty}^{\infty}(-1)^n
q^{2k\binom{n+1}{2}-n(i+j)}.
\end{eqnarray*}
We substitute $q\rightarrow q^{2k}$, 
$z\rightarrow -q^{-i-j}$ in \eqref{JTP} and
get
$$
\overline{\mathcal E}_{k,i}(1,q)=\frac{(-q)_\infty}{(q)_\infty}\sum_{j=0}^{2(k-i)}
(-1)^j (q^{i+j},q^{2k-i-j},q^{2k})_\infty.
$$

Finally we prove \eqref{gord1}.
$$
\overline{\mathcal E}_{k,i}(1/q,q)=\frac{(-1)_\infty}{(q)_\infty}\sum_{n=-\infty}^{\infty}
(-1)^n q^{2k\binom{n+1}{2}-ni}\frac{1+q^n}{2}.
$$
We substitute $q\rightarrow q^{2k}$, $z\rightarrow -q^{-i}$ and
$q\rightarrow q^{2k}$, $z\rightarrow -q^{-i+1}$ in \eqref{JTP} and
get
$$
\overline{\mathcal E}_{k,i}(1/q,q)=\frac{(-q)_\infty
}
{(q)_\infty}(q^{i},q^{2k-i},q^{2k};q^{2k})_\infty+(q^{i-1},q^{2k-i+1},q^{2k};q^{2k})_\infty.
$$

Now we give some combinatorial interpretation of Equations \eqref{gordon}, 
\eqref{gord} and \eqref{gord1}.

\subsection{2-modular diagrams}

We state in details the result for the Andrews' generalization of the  Gordon-G\"ollnitz identities
which correspond to  equation 
\eqref{gordon} of Corollary \ref{aftergf}. The coefficient of $q^n$
in this equation is the number of partitions of $n$ with parts 
not congruent to $2\mod 4$
and $0,\pm (2i-1)\mod 4k$. We make the change
of variable  $q\rightarrow q^2$ and $a\rightarrow 1/q$ in Theorem \ref{main}
and interpret it combinatorially in terms of 2-modular
diagrams defined in Section 2.
There exists an easy bijection $\phi$ between 2-modular diagrams of weight
$n$ with $j$ ones and overpartitions of $(n+j)/2$ with $j$ overlined 
parts. This bijection consists of erasing any 2 of the modular diagram and
changing any 1 to a marked corner. With this bijection in hand, the 
successive ranks (resp. Durfee dissection) of a 2-modular diagram $\mu$
 are the successive ranks (resp. Durfee dissection) of the corresponding
overpartition $\phi(\mu)$.

\begin{proposition}
All the following are equal:
\begin{itemize}
\item The number of partitions of $n$ with parts not congruent to $2\mod 4$
and $0,\pm (2i-1)\mod 4k$;
\item The number of partitions of $n$ of the form $(\lambda_1,\lambda_2,\ldots ,\lambda_s)$ with unrepeated odd parts, where
$\lambda_\ell - \lambda_{\ell+k-1} \ge 3$ if $\lambda_{\ell+k-1}$ is even
and $2$ otherwise; and $f_1+f_2<i$;
\item The number of 2-modular diagrams of $n$ whose successive ranks
lie in $[-i+2,2k-i-1]$;
\item The number of 2-modular diagrams of $n$ with $i-1$ successive
Durfee squares followed by $k-i$ successive Durfee rectangles, the first one being a generalized Durfee square/rectangle.
\item The number of paths that use four kinds of unitary steps with special
$(k,i)$-conditions where $n$ is twice the sum of the $x$-coordinates  of the peaks minus the number of South steps.
\end{itemize}
\end{proposition}

\noindent {\bf Remark.} 
The first two parts of the theorem are the Andrews' generalization of the Gordon-G\"ollnitz 
identities. The interpretation in terms of successive ranks and Durfee dissection 
is new to our knowledge.

\subsection{Superpartitions}

We give new partition theorems related to 
Gordon's theorems for overpartitions \cite{lo} which are the combinatorial 
interpretation of the case $i=k$ of equation (\ref{gord}) and $i=1$
of equation (\ref{gord1}).
We now interpret combinatorially equations (\ref{gord}) and (\ref{gord1}).
Superpartitions \cite{l} are overpartitions where the first
occurrence of a part can be overlined and the part $\overline{0}$
can appear. Let $\overline{B}_{k,i}(n)$ be the number of
overpartitions of $n$ of the form $(\lambda_1,\lambda_2,\ldots
,\lambda_s)$, where $\lambda_\ell - \lambda_{\ell+k-1} \ge 1$ if
$\lambda_{\ell+k-1}$ is overlined and $\lambda_\ell -
\lambda_{\ell+k-1} \ge 2$ otherwise and at most $i-1$ parts are
equal to $1$.
\begin{theorem}
For $1\le i\le k-1$,
the number of overpartitions counted by $\overline{B}_{k,i}(n)$ plus the
number of overpartitions counted by $\overline{B}_{k,i+1}(n)$ is equal to the
number of superpartitions where the non-overlined parts are not congruent to $0,\pm i$ modulo
$2k$.
\end{theorem}
\begin{proof}
Theorem \ref{main} tells us that 
$\overline{B}_{k,i}(n)+\overline{B}_{k,i+1}(n)$ is the coefficient of $q^n$
of $\overline{\mathcal E}_{k,i}(1,q)+\overline{\mathcal E}_{k,i+1}(1,q)$.

Thanks to equation \eqref{gord}, we know that 
$\overline{\mathcal E}_{k,i}(1,q)+\overline{\mathcal E}_{k,i+1}(1,q)=$
{\small \begin{eqnarray*}
%\overline{\mathcal E}_{k,i}(1,q)+\overline{\mathcal E}_{k,i+1}(1,q)&=&
&= & \frac{(-q)_\infty}{(q)_\infty}\left(\sum_{j=0}^{2(k-i)}
(-1)^j (q^{i+j},q^{2k-i-j},q^{2k};q^{2k})_\infty+\sum_{j=0}^{2(k-i-1)}
(-1)^j (q^{i+1+j},q^{2k-i-1-j},q^{2k};q^{2k})_\infty\right)\\
&=&
\frac{(-q)_\infty}{(q)_\infty}\left(\sum_{j=0}^{2(k-i)}
(-1)^j (q^{i+j},q^{2k-i-j},q^{2k};q^{2k})_\infty-\sum_{j=1}^{2(k-i)-1}
(-1)^j (q^{i+j},q^{2k-i-j},q^{2k};q^{2k})_\infty\right)\\
&=&\frac{(-1)_\infty}{(q)_\infty}(q^{i},q^{2k-i},q^{2k};q^{2k})_\infty.
\end{eqnarray*}}
The coefficient of $q^n$ in that last equation is the
number of superpartitions where the non-overlined parts are not congruent to $0,\pm i$ modulo
$2k$.
\end{proof}

\begin{theorem}
For $2\le i\le k-1$, the number of superpartitions of 
$n$ of the form $(\lambda_1,\lambda_2,\ldots ,\lambda_s)$, where
$\lambda_\ell - \lambda_{\ell+k-1} \ge 1$ if $\lambda_{j}$ is overlined and $\lambda_\ell - \lambda_{\ell+k-1} \ge 2$ otherwise and
where the number of non-overlined ones plus the number of $\overline{0}$ is at most $i-1$
is equal to the number of overpartitions of $n$  where the non-overlined parts are not congruent to $0,\pm i$ modulo $2k$
plus the number of overpartitions of $n$  where the non-overlined parts are not congruent to
$0,\pm (i-1)$ modulo $2k$.
\end{theorem}
\begin{proof}
We interpret combinatorially the coefficient of $q^n$ in equation  
\eqref{gord1}. This is the number of overpartitions of $n$  where the non-overlined parts are not congruent to $0,\pm i$ modulo $2k$
plus the number of overpartitions of $n$  where the non-overlined parts are not congruent to
$0,\pm (i-1)$ modulo $2k$. Note that this is the interpretation of 
$\overline{\mathcal E}_{k,i}(1/q,q)$. This implies that all the overlined parts
in Theorem \ref{main} are decreased by one and the result follows.
\end{proof}

\section{Conclusion}

We showed in this work how the combinatorial interpretation of the Andrews-Gordon
identities can be generalized to the case of overpartitions, when the combinatorial statistics 
(successive ranks, generalized Durfee square, length of the multiplicity sequence) are defined
properly. There exist other generalizations
of the Rogers-Ramanujan identities, see for example \cite{b1}. It was shown that the combinatorial
interpretation in terms of lattice paths can also be done for these identities \cite{agb,b,bu1,bu2}.
Our work can also be extended in that direction and the results are presented in \cite{clm}.
Recently Lovejoy and the second author have shown how to generalize the results
presented in this paper and in the paper \cite{clm} to overpartition pairs. This work appears in \cite{lm}.
Finally there exists
an extension of the concept of successive ranks for partitions due to Andrews, Baxter, Bressoud, Burge,
Forrester and Viennot \cite{abb} and our goal now is to extend that notion to overpartitions.\\

%\noindent{\bf Acknowledgments.} The authors  thank Jeremy Lovejoy for his help.

\bibliographystyle{alpha}

\begin{thebibliography}{A}
\bibitem{Ag-An-Br1}
A.K. Agarwal, G.E. Andrews, and D.M. Bressoud, The Bailey lattice,
\emph{J. Indian. Math. Soc.} {\bf 51} (1987), 57-73.
\bibitem{agb} A. K. Agarwal and D. M. Bressoud, Lattice paths and multiple basic
hypergeometric series. \textit{Pacific J. Math.} \textbf{135} (1989) 209--228.
\bibitem{ag} G. E. Andrews, The theory of partitions. Cambridge University Press, Cambridge, 1998.
\bibitem{ag2} G. E. Andrews, An analytic generalization of the Rogers-Ramanujan identities for odd moduli. \textit{Proc. Nat. Acad. Sci. USA} \textbf{71} (1974) 4082--4085.
\bibitem{ag3} G. E. Andrews, Sieves in the theory of partitions. \textit{Amer. J. Math.} \textbf{94} (1972) 1214--1230.
\bibitem{ag4} G. E. Andrews, A generalization of the G\"ollnitz-Gordon partition theorems. \textit{Proc. Amer. Math. Soc.} \textbf{8} (1967) 945--952.
\bibitem{abb} G. E. Andrews,  R. J. Baxter,  D. M. Bressoud, W. H. Burge, P. J. Forrester and G. Viennot,
Partitions with prescribed hook differences, \textit{European J. Combin.} (1987) \textbf{8} 341--350.
\bibitem{ab} G. E. Andrews and D. Bressoud, On The Burge correspondence between partitions and binary words. Number Theory (Winnipeg, Man., 1983). \textit{Rocky Mountain J. Math.} \textbf{15} (1985), no. 2, 225--233.
\bibitem{ab2}G. E. Andrews and D. M. Bressoud, Identities in Combinatorics III: further aspects of ordered set sorting. \textit{Discrete Math} \textbf{49} (1984) 222--236.
\bibitem{bep}
A. Berkovich and P. Paule, Lattice paths, $q$-multinomials and two
variants of the Andrews-Gordon identities, \emph{Ramanujan J.}
{\bf 5} (2001), 409--425.
\bibitem{at}
A.O.L. Atkin, A note on ranks and conjugacy of partitions, 
\textit{Quart. J. Math. Oxford Ser. (2)}  \textbf{17}  (1966) 335--338.
\bibitem{bp1} C. Bessenrodt and I. Pak, Partition congruences by involutions, \textit{ European J. Combin.} \textbf{25} (2004) 1139--1149.
%\bibitem{bp} C. Boulet and I. Pak, A combinatorial proof of the Rogers-Ramanujan and Schur identities. 
\textit{ J. Combin. Theory Ser. A},  \textbf{113}  (2006),  no. 6, 1019--1030.
\bibitem{br} F. Brenti, Determinants of Super-Schur Functions, Lattice Paths, and Dotted Plane Partitions
\textit{Adv. Math.}, \textbf{98} (1993), 27-64.
\bibitem{b} D. Bressoud, Lattice paths and the Rogers-Ramanujan identities. Number Theory, Madras 1987, 140--172, \textit{Lecture Notes in Math.} \textbf{1395}, Springer, Berlin, 1989.
\bibitem{b1} D. Bressoud, A generalization of the Rogers-Ramanujan identities for all moduli. \textit{J. Comb. Theory} \textbf{27} (1979) 64--68.
\bibitem{bu1} W. H. Burge, A correspondence between partitions related to generalizations of the Rogers-Ramanujan identities. \textit{Discrete Math.} \textbf{34} (1981), no. 1, 9--15.
\bibitem{bu2} W. H. Burge, A three-way correspondence between partitions. \textit{European J. Combin.} \textbf{3} (1982), no. 3, 195--213.
\bibitem{cl} S. Corteel, J. Lovejoy, Overpartitions. \textit{Trans. Amer. Math. Soc.} \textbf{356} (2004), no. 4, 1623--1635.
\bibitem{cl1} S. Corteel, J. Lovejoy, Frobenius partitions and the combinatorics of Ramanujan's $_1\psi_1$ summation, \textit{J. Comb. Theory Ser. A} \textbf{97} (2002), 177--183.
\bibitem{clm} S. Corteel, J. Lovejoy and O. Mallet, An extension to overpartitions of the Rogers-Ramanujan identities for even moduli,  preprint.
\bibitem{l}  P. Desrosiers, L. Lapointe, and P. Mathieu, 
 Jack superpolynomials, superpartition ordering and determinantal 
formulas.  \textit{Comm. Math. Phys.}  \textbf{233}  (2003),  no. 3, 383--402
\bibitem{fjm}J-F. Fortin, P. Jacob and P. Mathieu,  Jagged partitions.  
\textit{Ramanujan J.}  \textbf{10}  (2005),  no. 2, 215--235.
\bibitem{fm} J.-F. Fortin,  P. Jacob, and P.  Mathieu,
Generating function for $K$-restricted jagged partitions.
\textit{Electron. J. Comb.} \textbf{12}, No.1, Research paper R12, 17 p., (2005).
\bibitem{fm1} J.-F. Fortin,  P. Jacob, and P.  Mathieu,
Jagged partitions, \textit{Ramanujan J.} \textbf{10} (2005), no. 2, 215--235.
\bibitem{gr} G. Gaspar and M. Rahman, Basic hypergeometric series, Cambridge University
Press (1990).
\bibitem{go} B. Gordon, A combinatorial generalization of the Rogers-Ramanujan identities,
\textit{Amer. J. Math.} \textbf{83} (1961), 393--399 .
\bibitem{kw} S-J. Kang, J-H. Kwon,  Crystal bases of the Fock space representations and string functions.  
\textit{J. Algebra}  \textbf{280}  (2004),  no. 1, 313--349.
\bibitem{lo} J. Lovejoy, Gordon's theorem for overpartitions. \textit{J. Combin. Theory Ser. A} \textbf{103} (2003), no. 2, 393--401.
\bibitem{lo1} J. Lovejoy, Overpartition theorems of the Rogers-Ramanujan type
\textit{J. London Math. Soc.} \textbf{69} (2004), 562--574.
\bibitem{lo2} J. Lovejoy, Rank and conjugation for the Frobenius representation of an overpartition. \textit{Ann. Combin.} \textbf{9} (2005) 321--334.
\bibitem{lm} J. Lovejoy and O. Mallet, Overpartition pairs and two families of basic hypergeometric series, preprint.
\bibitem{mm} P. A. MacMahon, Combinatory Analysis, vol. 2, Cambridge Univ. Press, Cambridge, 1916.
\bibitem{ma} K. Mahlburg, The overpartition function modulo small powers of 2, \textit{Discrete Math.} 
\textbf{286} (2004), no. 3, 263--267.
\bibitem{m} O. Mallet, Rangs successifs, chemins et identit\'es de type Rogers-Ramanujan, \textit{Master's thesis, Paris 6} (2005).
\bibitem{MAL} O. Mallet, \textit{PhD thesis}, Paris 7, in preparation.
\bibitem{ro} L. J. Rogers, Second memoir on the expansion of certain infinite products, \textit{Proc. London Math. Soc.} \textbf{25} (1894), 318--343.
\bibitem{Sl1}
L.J. Slater, A new proof of Rogers's transformations of infinite
series \emph{Proc. London Math. Soc.} {\bf 53} (1951), 460-475.
\bibitem{y} A. J. Yee, Combinatorial proofs of Ramanujan's $_1\psi_1$ summation and the $q$-Gauss summation, \textit{J. Combin. Theory Ser. A}, \textbf{105} (2004), 63--77.
\end{thebibliography}

\end{document}